\pdfoutput=1
\RequirePackage{ifpdf}
\ifpdf 
\documentclass[pdftex]{sigma}
\else
\documentclass{sigma}
\fi

\usepackage{tikz,tikz-cd}
\usepackage[all]{xy}

\numberwithin{equation}{section}
\newtheorem{Theorem}{Theorem}[section]

\newtheorem{Lemma}[Theorem]{Lemma}
\newtheorem{Proposition}[Theorem]{Proposition}
 { \theoremstyle{definition}
\newtheorem{Definition}[Theorem]{Definition}

\newtheorem{Example}[Theorem]{Example}
\newtheorem{Remark}[Theorem]{Remark}
\newtheorem{notation}[Theorem]{Notation}
\newtheorem{construction}[Theorem]{Construction}
}

\newcommand{\vv}{\mathrm{vir}}

\newcommand{\E}{\mathfrak E}

\newcommand{\OO}{\mathcal{O}}

\newcommand{\rk}{\operatorname{rk}}

\newcommand{\Spec}{\operatorname{Spec}}
\newcommand{\Proj}{\operatorname{Proj}}

\newcommand{\PP}{\mathbb P}

\newcommand{\MM}[4]{\overline{\mathcal M}_{#1,#2}(#3,#4)}

\begin{document}

\allowdisplaybreaks

\newcommand{\arXivNumber}{1804.06048}

\renewcommand{\thefootnote}{}

\renewcommand{\PaperNumber}{026}

\FirstPageHeading

\ShortArticleName{Virtual Classes for the Working Mathematician}

\ArticleName{Virtual Classes for the Working Mathematician\footnote{This paper is a~contribution to the Special Issue on Integrability, Geometry, Moduli in honor of Motohico Mulase for his 65th birthday. The full collection is available at \href{https://www.emis.de/journals/SIGMA/Mulase.html}{https://www.emis.de/journals/SIGMA/Mulase.html}}}

\Author{Luca BATTISTELLA~$^\dag$, Francesca CAROCCI~$^\ddag$ and Cristina MANOLACHE~$^\S$}

\AuthorNameForHeading{L.~Battistella, F.~Carocci and C.~Manolache}

\Address{$^\dag$~Mathematikon, Im Neuenheimer Feld 205, 69120 Heidelberg, Germany}
\EmailD{\href{mailto:lbattistella@mathi.uni-heidelberg.de}{lbattistella@mathi.uni-heidelberg.de}}

\Address{$^\ddag$~JCMB, Peter Guthrie Tait Road, Edinburgh EH9 3FD, UK}
\EmailD{\href{mailto:francesca.carocci@ed.ac.uk}{francesca.carocci@ed.ac.uk}}

\Address{$^\S$~Hicks building, 226 Hounsfield Rd, Broomhall, Sheffield S3 7RH, UK}
\EmailD{\href{mailto:c.manolache@sheffield.ac.uk}{c.manolache@sheffield.ac.uk}}

\ArticleDates{Received October 03, 2019, in final form March 25, 2020; Published online April 09, 2020}

\Abstract{This note is intended to be a friendly introduction to virtual classes. We review virtual classes and we give a number of properties and applications. We also include a new virtual push-forward theorem and many computations of virtual classes in simple examples.}


\Classification{14C17; 14N35; 14D23}

\renewcommand{\thefootnote}{\arabic{footnote}}
\setcounter{footnote}{0}

\section{Introduction}
Virtual classes have been introduced by Li--Tian \cite{litian} and Behrend--Fantechi \cite{BerFan} and are a key tool in enumerative geometry. However, virtual classes are usually considered to be difficult, possibly because there are very few explicit computations. The goal of this note is to show that virtual classes of well understood schemes are usually computable. The main reason for the lack of computations is not the definition of virtual classes -- which may indeed look unsuited for computations -- but the poor understanding of spaces on which we are usually interested to compute virtual classes. The main examples of spaces with virtual classes are moduli spaces of curves on a given variety such as moduli spaces of stable maps, Hilbert schemes, etc. All these spaces have the same major drawback: we typically know very little about their geometry. We expect that if we had a space with known cohomology and equations of this space inside a smooth ambient space, then the computation of its virtual class would be possible (see Section~\ref{algorithms}).

We now briefly describe the content of this paper. One of the ingredients in the construction of virtual classes is the normal cone. In Section~\ref{normalcones} we give examples of explicit computations of normal cones and we list properties of normal cones and their Segre classes. The second ingredient is a choice of a perfect obstruction theory. We define obstruction theories first in Section~\ref{gysinmorphism} and then more generally in Section~\ref{gluingobs}. In Section~\ref{gluingobs} an obstruction theory $E_{X/V}$ on a~scheme or DM stack $X$ is given by the following data:
\begin{itemize}\itemsep=0pt
\item a global embedding $X\hookrightarrow V$, with $V$ a smooth scheme or DM stack,
\item a covering of $V$ by open sets $V_i$, and on each $U_i=X\cap V_i$ a set of generators for the ideal~$I_{U_i/V_i}$,
\item gluing data for the sets of generators.
\end{itemize}
We compute examples of virtual classes associated to such data. In the end of Section~\ref{gluingobs} we review the more general construction in~\cite{BerFan}. The main advantage of the construction in~\cite{BerFan} is that one does not need a global embedding of $X$ in a smooth scheme/ stack. In most examples in this note, we will assume such an embedding exists. This is not a very restrictive assumption for our purposes; for example moduli spaces of stable maps admit such embeddings \cite[Appendix~A]{grpa}. In Section~\ref{deformations} we review a variant of the deformation to the normal cones in \cite{KKP} and we use this construction to compute virtual classes in certain examples. In Section~\ref{properties} we list properties of virtual classes and we give two computational methods. Section~\ref{properties} includes a new virtual push forward theorem for $X\to Y$, with $X$ a component of $Y$ (Theorem \ref{vpfnew}). In Section~\ref{applications} we give examples of moduli spaces of stable maps. We discuss the hyperplane property for Gromov--Witten invariants in genus 0 and 1. For genus 1 stable maps we suggest a way to split the virtual class on the components of the moduli space in a simple example.

\subsection*{Old and new}

For Chow groups, normal cones, Gysin maps and their properties the main reference is Fulton's intersection theory book \cite{Fulton}. We borrow most of our examples from it. All Chow groups will be Chow groups with rational coefficients. For Chow groups of Deligne--Mumford stacks we refer to \cite{vis}. The stacky version of all the above is explained in \cite{BerFan}. The definition of virtual classes which is closest to our approach is the one of Siebert \cite{siebert}. We have tried to give an overview of some of the main constructions and techniques in virtual intersection theory, but no attempt has been made to give an exhaustive survey of this topic. There are many fundamental results which we do not treat here. Among these, we mention the virtual localization theorem~\cite{grpa} of Graber and Pandharipande, which is by far the most used method to compute virtual classes.

While this note is intended to be expository, there are some results with a certain degree of originality. The virtual push-forward Theorem~\ref{vpfnew} and some ideas in Section~\ref{algorithms} are new. Section \ref{applications} also contains some new examples and observations.

\section{Motivation}

\subsection{Enumerative geometry} Suppose we want to count a certain type of smooth curves in a smooth space $X$. The counting problem is usually solved in two steps: first we construct a proper parameter space for our curves in $X$, and then we do intersection theory on this space.

There are by now many proper moduli spaces of curves on a given variety. Grassmannians have been the first moduli spaces used to solve enumerative problems. Modern examples of moduli spaces of curves on a given variety include Hilbert schemes and moduli spaces of stable maps. For a more complete list see \cite{rt}.

Such moduli spaces are often singular, with many different irreducible components, possibly of different dimensions. Geometric conditions can often be translated into cohomology classes, but we need them to act on a fundamental class in order to get numerical invariants. The solution to this problem is to construct \emph{virtual fundamental classes} -- these mimic the role of the fundamental cycle in the smooth/pure dimension case.

After having defined intersection theory on these spaces, we see that we sometimes get wrong answers to enumerative problems. This is because we also count \emph{degenerate contributions} from curves which are not smooth, or have different numerical invariants. These degenerate contributions come with different flavours: they could lie on the boundary of our moduli space, or they could form entire ``boundary components''.

Virtual classes and degenerate contributions are the main concern of these notes. Let us look at two elementary examples.

\subsection{Conics tangent to five lines}\label{conics five lines} Let us determine the number of conics tangent to five general lines in the plane. This is an example in which the moduli space is irreducible and smooth, but the intersection is non-transverse and it gives us the wrong answer.

A plane conic is given by a quadratic polynomial
\[\sum_{0\leq i\leq j\leq2}a_{i,j}x_ix_j\]
up to scalar multiplication. The parameter space of conics in $\mathbb P^2$ is thus identified with $\mathbb P^5$ with homogeneous coordinates $[a_{00}:a_{10}:\dots :a_{22}]$.

To express the tangency condition, consider the line $L=\{x_0=0\}$. A conic is tangent to $L$ if and only if the polynomial $a_{11}x_1^2+a_{12}x_1x_2+a_{22}x_2^2$ has a double root, i.e.,
\[
\Delta_{L}=4a_{11}a_{22}-a_{12}^2=0.
\]
This gives a quadric $Q$ in $\mathbb P^5$. We might expect that the number of conics tangent to five lines in general position is given by intersecting the corresponding five quadrics. If $H$ denotes the hyperplane class in $A_*\big(\mathbb P^5\big)$, then the intersection of five general quadrics is \begin{gather*}(2[H])^5=2^5=32.\end{gather*} This is wrong: the correct number of conics tangent to five given lines is 1 (dually, it is the number of conics through five points). The reason why this count gives the wrong answer is that the intersection $\bigcap_{i=1}^5 Q_i$ is not transverse: it contains the positive-dimensional locus of double lines $Z\simeq \mathbb P^{2*}$. The contribution of~$Z$ to the intersection product $(2[H])^5$ can be computed via residual intersections (see \cite[Chapter~9]{Fulton}). We compute this contribution in Example~\ref{conics lines revised}.

\subsection{A line and a plane}
Let us now look at a virtual class on a space with components of different dimensions. The algebraic construction will be done in Section \ref{gysinmorphism}. This example is modelled on a certain moduli space of stable maps (see Section \ref{applications}).

Consider $X=V(xz,yz)\subset\PP^3$, which is the union of the plane $H=(z=0)$ and the line $L=(x=y=0)$. Note that $X$ is the zero locus of two non transverse sections of $\OO_{\PP^3}(2)$. If the two sections were transverse, their intersection would have been a curve of degree~4 in $\PP^3$. We would like to define a virtual class in $A_*(X)$ which is related to the transverse intersection. Such a class on $X$ can be constructed in the following way. Deforming the equations which define $X$ to $\big(xz-\epsilon y^2,yz-2\epsilon x^2\big)$ we get the union $L\cup C_{\epsilon}$, where $C_{\epsilon}$ is a twisted cubic in $\mathbb P^3$. Taking the limit for $\epsilon$ going to zero, we get the union of the class of $L$ and the class of a degree~3 plane curve. Since virtual classes are supposed to be deformation invariants, the class of this limit in~$A_*(X)$ is a good candidate for the virtual class of~$X$. We will denote it by $[X]^{\vv}$.

However, deformations are usually hard to keep track of and they would make our virtual class unsuited for computations in most cases of interest. In the following we will not construct a class on a deformation of $X$. Instead, we will use Fulton--Macpherson intersection theory to define a class on our reducible space.

\begin{Remark} Let $X=V(xz,yz)\subset\PP^3$ as above. If $i:X\to\PP^3$ denotes the embedding, then we can define $i_*[X]^{\vv}$ by
\[
i_*[X]^{\vv}:=c_2(\mathcal{O}(2)\oplus \mathcal{O}(2))\cap\big[\PP^3\big].
\]
The push-forward of $[X]^{\vv}$ is easy to compute and if we are interested in degree zero intersections of $[X]^{\vv}$ with classes which are pulled back from classes on $\PP^3$, we can work with $i_*[X]^{\vv}$ instead of $[X]^{\vv}$.

We do not find this satisfactory. One reason is that we may be interested in separating $[X]^{\vv}$ on the components of $X$. This information is lost after pushing forward to $\PP^3$.
\end{Remark}

\section{Normal cones}\label{normalcones}

The main ingredient in Fulton's intersection theory is the normal cone of a closed embedding. We review normal cones and their properties and we give many examples. The main reference for this section is \cite[Chapters~4 and~5]{Fulton}. We assume familiarity with Chow groups (see \cite[Chapters~1 and~2]{Fulton}).

\begin{Definition}
Let $X\hookrightarrow Y$ be a closed embedding and $\mathcal I=\mathcal{I}_{X/Y}$ the ideal sheaf of $X$ in $Y$. We define the
\emph{normal cone} of $X$ in $Y$ as
\[C_{X/Y}=\underline{\operatorname{Spec}}_X\bigg(\bigoplus_{n\geq 0}\mathcal{I}^n/\mathcal{I}^{n+1}\bigg)\]
and the \emph{normal sheaf} of $X$ in $Y$ as
\[N_{X/Y}=\underline{\operatorname{Spec}}_X\big({\mathcal S}{\rm ym}^\bullet \big(\mathcal{I}/\mathcal{I}^2\big)\big).\]
\end{Definition}

\begin{Example}
If $X\hookrightarrow Y$ is a regular embedding, then $C_{X/Y}=N_{X/Y}$ and it is a vector bundle. Furthermore, for $X$ and $Y$ smooth, $N_{X/Y}$ fits in the exact sequence of vector bundles:
\[0\to T_X\to T_{Y}\rvert_X\to N_{X/Y}\to 0.\]
\end{Example}

\begin{Example}\label{main example cone}
Consider $X=V(xz,yz)\subseteq\PP^3_{[x:y:z:w]}$. In the following we compute $C_{X/\PP^3}$. In the affine chart $w\neq 0$ the coordinate ring of $X$ is $R=\mathbb C[x,y,z]/(xz,yz)$. In this chart, the normal cone is
\[
C_{X/\PP^3}\rvert_{U_w}=\operatorname{Spec}\left(R[A,B]/(yA-xB)\right).
\]
To see this, consider the surjective morphism
 \begin{align*}
 R[A,B] & \rightarrow \bigoplus_{n\geq 0}{I^n_{X/\PP^3}/I^{n+1}_{X/\PP^3}},\\
 A &\rightarrow xz,\\
 B &\rightarrow yz.
 \end{align*}
Then, $C_{X/\PP^3}\rvert_{U_w}\subset X\rvert_{U_w}\times \mathbb{A}^2$ is given by the kernel of this map, which is generated by $(yA-xB)$. Local descriptions of $C_{X/\PP^3}$ glue to give an embedding
 \[
 C_{X/\PP^3}\hookrightarrow\mathcal O_X(2)\oplus\mathcal O_X(2).
 \]

The cone $C_{X/\PP^3}\to X$ has two components as follows.
\begin{itemize}\itemsep=0pt
\item Over the line L, the cone is defined by the ring $R[A,B]$. This shows that over the $z$-axis, $C_{X/\PP^3}\to X$ is isomorphic to the rank 2 vector bundle $\mathcal{O}_L(2)\oplus\mathcal O_L(2)$.
\item Over the plane $H$ we have two different situations: if at least one between $x$ and $y$ is non zero, say $y$, then the defining ring of the cone is $R[B]$, i.e. $C_{X/\PP^3}\to X$ is a line bundle; over the origin the fibre has rank 2.
\end{itemize}
We have that $C_{L/\PP^3}=\mathcal{O}_L(1)\oplus\mathcal O_L(1)$ and the natural map $C_{L/\PP^3}\to C_{X/\PP^3}$ collapses the fiber over the origin $z=0$. This shows that
\[C_{X/\PP^3}\neq C_{L/\PP^3}\cup C_{H/\PP^3}.\]

\end{Example}

Let us now list some general properties of normal cones.

\subsection{Components of the cone} The main reference is \cite[Section~4.1]{Fulton}.
\begin{enumerate}\itemsep=0pt
\item If $Y$ has irreducible components $Y_1,\dots, Y_N$, then
 \[C_{X/Y}=\bigcup_i^N C_{X_i/Y_i},\qquad X_i=Y_i\times_Y X.\]
 \item The number of components of $C_{X/Y}$ is always bigger or equal than the number of components of $X$, and it may be strictly bigger.
\end{enumerate}
 \begin{Example}Consider the inclusion of the singular point in the nodal cubic $X=pt\hookrightarrow Y=V\big(y^2+x^2(x-1)\big)$. Then $X$ has only one component and one can easily compute that $C_{X/Y}=V\big(y^2-x^2\big)$ which has two components.
 \end{Example}

 \begin{Example}\label{dorigin}
Let us compute the normal cone of $X=V\big(x^2,xy\big)\subset\mathbb{A}^2$. Note that $X$ is the $y$-axis with an embedded point in the $x$ direction at the origin.

There is a graded surjective morphism
\begin{align*}
\mathbb{C}[x,y]/\big(x^2,xy\big)[A,B]&\to\bigoplus_{n\geq 0}{ I^n_{X/\mathbb{A}^2}/ I^{n+1}_{X/\mathbb{A}^2} },\\
A &\to x^2,\\
B &\to xy.
\end{align*}
To understand the ideal of $C_{X/\mathbb{A}^2}\subseteq X\times\mathbb{A}^2$ we need to compute $J$ -- the kernel of this morphism.
We see that $J$ is generated by $yA-xB$ and thus
\[
C_{X/\mathbb{A}^2} =\Spec\left( \frac{\mathbb{C}[x,y]}{\big(x^2,xy\big)}[A,B]/(yA-xB) \right).
\]
Since the only relation we have is linear in $A$ and $B$, the normal cone coincides with the normal sheaf.

Let us now look at the components of $C_{X/\mathbb{A}^2}$. Over the generic point of the $y$-axis, the fibre of $C_{X/\mathbb{A}^2}\to X$ is $\mathbb{A}^1$, so $C_{X/\mathbb{A}^2}$ restricted to the punctured line is a line bundle $L$. The fibre over the origin is isomorphic to $\mathbb{A}^2$. In conclusion, $C_{X/\mathbb{A}^2}$ has two irreducible components, one supported on the origin and another one which is the closure of $L$ in $X\times \mathbb{A}^2$. See Example \ref{doublepoint} for more details.
\end{Example}

\subsection{Pull-backs}
Given a cartesian diagram of schemes
\[\xymatrix{X'\ar[r]\ar[d]^f &Y'\ar[d]\\
			X\ar[r] & Y}\]
there is a closed embedding $C_{X'/Y'}\hookrightarrow f^* C_{X/Y}.$

\subsection{Deformation to the normal cone}\label{sec:defnormal} The main reference is \cite[Chapter~5]{Fulton}. Let $X$ be a closed subscheme of $Y$ and let $C_{X/Y}$ be the normal cone of the embedding. There exists a flat family $\pi\colon M^0_{X/Y}\to \PP^1$ with generic fiber $Y$ and special fiber $C_{X/Y}$ together with a closed embedding
\[X\times \PP^1\hookrightarrow M^0_{X/Y}\]
 commuting with the maps to $\PP^1$ such that
 \begin{itemize}\itemsep=0pt
 \item over $\PP^1\backslash \{0\}$ we get the given embedding of $X$ in $Y$, and
 \item the fiber over $0$ is the embedding of $X$ into its normal cone as the zero section.
 \end{itemize}
To construct this deformation, we first define
\[M_{X/Y}=\operatorname{Bl}_{X\times \{0 \}} \big(Y\times \PP^1\big).\]
One can check that the fiber of $M_{X/Y}$ over $t=0$ is $\PP_X(C_{X/Y}\oplus \mathcal O_X)\cup\operatorname{Bl}_X Y$ with the two components intersecting in $\PP(C_{X/Y})$. We then define $M^0_{X/Y}$ to be the complement of $\operatorname{Bl}_X Y$ in~$M_{X/Y}$.
\begin{Remark}
From the previous construction, we see that if $Y$ has pure dimension~$d$, then~$C_{X/Y}$ has pure dimension~$d$.
\end{Remark}
\subsection{Segre classes}
The main reference is \cite[Chapter~4]{Fulton}. To every cone we can associate a cycle called its Segre class. Methods to compute Segre classes are reviewed in Section~\ref{algorithms}.
\begin{Definition}
Let $X\hookrightarrow Y$ be an embedding and $C_{X/Y}$ the normal cone of~$X$ in~$Y$. Let $\mathcal O(1)$ be the tautological line bundle on
\[\PP(C_{X/Y}\oplus \mathcal O_X)\stackrel{q}{\longrightarrow}X.
\]
We define the \emph{Segre class} $s(X,Y)$ of $C_{X/Y}$ to be the cycle
\[
s(X,Y)=q_*\bigg(\sum_{i\geq 0}c_1(\mathcal O(1))^i\cap[\PP(C_{X/Y}\oplus \mathcal O_X)]\bigg).
\]
For a more general cone $q\colon C\to X$, one can define its Segre class by the analogous formula $s(C)=q_*\big(\sum_{i\geq 0}c_1(\mathcal O(1))^i\cap[\PP(C\oplus \mathcal O_X)]\big)$. See \cite[Chapter~4]{Fulton} for more details.
\end{Definition}
\begin{Example}\label{ex:s=cinv} If $C_{X/Y}=N_{X/Y}$ is a vector bundle, then
\[s(X,Y)=c(C_{X/Y})^{-1}.\]
Here $c$ denotes the total Chern class. (See \cite[Proposition~4.1]{Fulton}.)
\end{Example}
\begin{Example}\label{conics lines revised}
Let us look again at the example in Section \ref{conics five lines} and see that there is one conic tangent to five general lines in a plane. We use \cite[Proposition~9.1.1]{Fulton} and straightforward computations of Segre classes of vector bundles. Let $\nu_2\colon\mathbb P^2\to\mathbb P^5$ be the Veronese embedding of degree $2$ and let $Z$ denote the image of~$\mathbb P^2$. We can check that $Z$ parametrises the space of double lines inside the space of conics. Both $Z$ and $\PP^5$ are smooth and the Chern classes of their tangent bundles are easily computed from the Euler sequence. Let $L$ be the class of a line in $Z$, and $N_i$ the restriction of the normal bundle $N_{Q_i/\mathbb P^5}\simeq\mathcal O_{Q_i}(2)$ to $Z$. We have $c(N_i)=(1+4L)$. The contribution of $Z$ to $\deg\big(\bigcap_{i=1}^5 Q_i\big)$ is then computed in terms of the Segre class of the cone~$C_{Z/\PP^5}$ and the obstruction bundle $\oplus_{i=1}^5N_{Q_i/\mathbb P^5}$. Using the normal sequence of $Z$ in~$\PP^5$, we have that the Segre class of the cone $C_{Z/\PP^5}$ equals $c(T_{\mathbb P^5|Z})^{-1} c(T_Z)\cap [Z]$. We get that the contribution of $Z$ is given by
\begin{align*}
 (Q_1\cdots Q_5\cdot\mathbb P^5)^Z& =\left\{ \left(\prod c(N_i)\right) c(T_{\mathbb P^5|Z})^{-1} c(T_Z)\cap [Z]\right\}_0\\& =\deg\big((1+4L)^5(1+2L)^{-6}(1+L)^3\big[\mathbb P^2\big]\big) =31,
\end{align*}
where the notation $\{a\}_0$ means the degree zero part of the class $a\in A_*(Z)$. With this we see that the plane of double lines contributes~$31$ to the intersection $\cap_{i=1}^5Q_i$ and the residual scheme is just one reduced point which represents the nondegenerate conic tangent to the five given general lines.
\end{Example}

\begin{Proposition}[{\cite[Examples~4.1.6, 10.1.10]{Fulton}}]\label{Segre constant}
Let $\mathcal{C}$ be a cone over $X\times \mathbb{A}^1$ flat over $\mathbb{A}^1$ and let $C_t$ be the restriction of $\mathcal C$ to $X\times\{t\}$. Then
\[
s(C_0)=s(C_1) \in A_*(X).
\]
\end{Proposition}

{\bf Computing Segre classes by blowing up.} We recall some useful formulae from~\cite{Fulton}.

\begin{Proposition}[{\cite[Corollary 4.2.2]{Fulton}}]\label{prop2}
Let $X$ be a proper closed subscheme of a variety~$Y$, $ \tilde{Y}=\operatorname{Bl}_X Y$ and $\eta\colon\tilde{X}\to X$ the exceptional divisor. Then
\begin{gather}\label{prop2-eq}s(X,Y)=\sum_{k\geq 1}{(-1)^{k-1}\eta_*[\tilde{X}]^k}.\end{gather}
\end{Proposition}

\begin{Lemma}[{\cite[Example~8.3.9]{Fulton}}]\label{blup_formulas}
Let $X\subseteq Y$ be non-singular varieties, and consider $\tilde{Y}=\operatorname{Bl}_X Y$ with exceptional divisor $\tilde{X}$:
\[
\begin{tikzcd}
\tilde{X}\ar[r,"j"]\ar[d,"\eta"] & \tilde{Y}\ar[d,"f"] \\
X\ar[r,"i"] & Y
\end{tikzcd}
\]
Then the ring structure on $A_*\big(\tilde{Y}\big)$ is given by
\begin{enumerate}\itemsep=0pt
\item[$(i)$] $f^*y\cdot f^* y'=f^*(y\cdot y')$,
\item[$(ii)$] $j_*(\tilde{x})\cdot j_*(\tilde{x}')=j_*\big(c_1\big(j^*\mathcal O_{\tilde{Y}}\big(\tilde{X}\big)\big)\cdot \tilde{x}\cdot \tilde{x}'\big)$,
\item[$(iii)$] $f^*y\cdot j_*(\tilde{x})=j_*((\eta^* i^* y)\cdot\tilde{x})$.
\end{enumerate}			
\end{Lemma}

\subsection{Failure of exact sequences}	
Given a composition of regular embeddings \begin{gather*}X\hookrightarrow Y\hookrightarrow Z\end{gather*} there exists a short exact sequence of vector bundles
\[0\to N_{X/Y}\to N_{X/Z}\to N_{Y/Z}\rvert_X\to 0.\]
When the embeddings are not regular, we still have induced maps of cones
\begin{equation}\label{eq:exactcones}
C_{X/Y}\to C_{X/Z}\to C_{Y/Z}\rvert_X.\end{equation}
Intuitively, a sequence of cones should be called exact if locally we have $C_{X/Z}\simeq C_{X/Y}\times_X C_{Y/Z}\rvert_X$. See \cite[Example~4.1.6]{Fulton} for a more restrictive definition of exactness. Sequence~\eqref{eq:exactcones} is not exact in general.

\begin{Example}Exact sequences of cones
\begin{enumerate}\itemsep=0pt
\item[(a)] Consider $X\overset{j}{\hookrightarrow} Y \overset{i}{\hookrightarrow} Z=Y\times\mathbb A^n$, then $C_{X/Z}\simeq C_{X/Y}\times\mathbb A^n$.
\item[(b)] Let $X\overset{j}{\hookrightarrow} Y=X\times\mathbb A^n \overset{i}{\hookrightarrow} Z$. Then, the sequence of cones
\[0\to\mathbb A^n_X\to C_{X/Z}\to j^*C_{Y/Z}\to 0
\] is not necessarily exact.
\end{enumerate}
\end{Example}

\begin{proof}
(a) Locally we may assume
\begin{gather*}Y=\Spec R,\qquad X=\Spec R/I,\qquad Z=\Spec R[t_1,\ldots,t_n].
\end{gather*} Then the ideal of $X$ in $Z$ is $J=I+(t_1,\ldots,t_n)$, and intuitively there are no mixed relations between the generators of $J$ and the $t_i$'s. Indeed, there is a surjective map
\[
 \bigoplus_{i=0}^k I^i/I^{i+1}\otimes_{R/I} (t_1,\ldots,t_n)^{k-i}\to J^k/J^{k+1}.
\]
To show injectivity, pick a homogeneous polynomial of degree $k$ in the generators of $I$ and $t_i$. Suppose its image lies in $J^{k+1}$and the degree in $t_i$ is $k-j$. Then, the corresponding coefficient must then lie in $I^{j+1}$ which shows that the polynomial is zero on the left-hand side too.

(b) A counterexample is given by $Z$ the cone over a quadric, $Y$ a line through the vertex, and $X$ the vertex. One sees that the sequence
\[\mathbb C[v_0,v_2]\hookrightarrow \mathbb C[v_0,v_1,v_2]/\big(v_2^2-v_0v_1\big)\twoheadrightarrow\mathbb C[v_1]\]
does not split because
\[\mathbb C[v_0,v_2]\otimes_{\mathbb C}\mathbb C[v_1]\ncong \mathbb C[v_0,v_1,v_2]/\big(v_2^2-v_0v_1\big).\]
See Example~\ref{segrecones}(b) below for a proof with Segre classes.
\end{proof}

\begin{Remark} Case (a) admits the following generalisation (see \cite[Theorem~6.5]{Fulton}). For
 \[
 X\overset{j}{\hookrightarrow} Y \overset{i}{\hookrightarrow} E,
 \]
 with $E$ a vector bundle on $Y$ and $i$ the zero section, we have
 \[
 C_{X/Z}\simeq C_{X/Y}\times_X j^*E.
 \]
 This does not generalise to $E$ any cone on~$Y$.
\end{Remark}

\begin{Example}\label{segrecones}
 This is expanding on \cite[Example~4.2.8]{Fulton}.
\begin{enumerate}\itemsep=0pt
\item[(a)] Consider $X\overset{j}{\hookrightarrow} Y \overset{i}{\hookrightarrow} Z$ with $i$ a regular embedding; then not necessarily $s(X,Y)=s(X,Z)\cap c(N_{Y/Z})$.
\item[(b)] Consider $X\overset{j}{\hookrightarrow} Y \overset{i}{\hookrightarrow} Z$ with $j$ a regular embedding of codimension~1, i.e., $X\subseteq Y$ a~Cartier divisor; it is not always true that \begin{gather*}j^*s(Y,Z)=s(X,Z)\cap c(N_{X/Y}).\end{gather*}
\end{enumerate}
\end{Example}
\begin{proof} (a) Take $Z=\mathbb P^2$, $Y$ a nodal cubic, e.g., the one given by \begin{gather*}\big\{f=y^2z-x^3-x^2z=0\big\}\end{gather*} and $X$ the node, which is the origin of the affine chart $\mathbb A^2_{x,y}$. Then $Y\subseteq Z$ is a~Cartier divisor with normal bundle~$\mathcal O_Y(3)$. $C_{X/Z}$ is a vector bundle over $X$ and $s(X,Z)=[\underline 0].$ On the other hand, $s(X,Y)=2[\underline 0]$. This can be seen directly from the local isomorphism
\[C_{X/Y}\simeq \Spec\mathbb C[t_1,t_2]/\big(t_1^2-t_2^2\big).
\]
The cone $\Spec\mathbb C[t_1,t_2]/\big(t_1^2-t_2^2\big)$ has two irreducible components (the tangent directions to the two branches of the curve at the node).

(b) Let $Z$ be the affine cone given by \begin{gather*}\Spec \mathbb C[x,y,z]/\big(x^2+y^2-z^2\big).\end{gather*} Let $X$ be the vertex of the cone and $Y=\{x=y-z=0\}$ a line through $X$. Then the normal cone $C_{Y/Z}$ is given by \begin{gather*}\mathbb C[x,y,z][w_0,w_1]/(x,y-z,xw_0+(y+z)w_1).\end{gather*} Let us compute the Segre class separately of $C_{Y/Z}$. When restricting to the point, only the degree~1 part of $s(Y,Z)$ matters, hence we may discard the component supported on $\{x=y=0\}$. Since the remaining component is a vector bundle on $Y$, we get $j^*s(Y,Z)=[X]$.

Let us compute $s(X,Z)$ by means of formula~\eqref{prop2-eq}. Consider $\tilde Z=\operatorname{Bl}_X Z$ and $\tilde X$ the exceptional divisor. $Z$ has an $A_1$-singularity at $X$ and $\tilde Z$ is the minimal resolution. This shows that~$\tilde X$ is a $(-2)$-curve and thus
\[
 s(X,Z)=\pi_*\big([\tilde X]-\tilde X\cdot[\tilde X]\big)=0-(-2[X]).
\]
We thus get $j^*s(Y,Z)\neq s(X,Z)\cap c(N_{X/Y})$.
\end{proof}

\section{Virtual pull-backs}\label{gysinmorphism}
In order to define the virtual fundamental class, we need two more ingredients: the Gysin map and a vector bundle containing our cone. We define virtual classes in a simple context and we give examples. The main references are \cite{vir-pull} for the definitions, and \cite[Chapters 5 and 6]{Fulton} for the methods (in the schematic setting).

\begin{Proposition}[{{\cite[Theorem 3.3]{Fulton}}}] Let $E$ be a vector bundle of rank $r$ on $X$, then
the flat pull-back $\pi^*\colon A_{k-r}(X)\to A_k(E)$ is an isomorphism for any $k$.
\end{Proposition}

\begin{notation} \label{inverse} Let $0_E^!\colon A_k(E)\to A_{k-r}(X)$ denote the inverse of the flat pull-back
\[\pi^*\colon \ A_{k-r}(X)\to A_k(E).
\]
\end{notation}
\begin{construction}\label{naive}
Suppose that we have the following data:
\begin{enumerate}\itemsep=0pt
\item[1)] an embedding $X\hookrightarrow Y$, and
\item[2)] a vector bundle $E_{X/Y}$ of rank $r$ with an embedding of cones over $X$
 \begin{equation*}C_{X/Y}\hookrightarrow E_{X/Y}.
 \end{equation*}
\end{enumerate}
 Then, we have morphisms:
\begin{align}\label{virtpb}
A_*(Y)\xrightarrow{\sigma} A_*(C_{X/Y})\xrightarrow{i_*} A_*(E_{X/Y})\xrightarrow{0_{E_{X/Y}}^!} A_{*-r}(X),
\end{align}
where the first arrow is defined as $\sigma[V]=[C_{V\cap X} V]$, the second one is the proper push-forward, and the third one is the morphism in Notation~\ref{inverse}. We denote the composition in \eqref{virtpb} by $f^!_{E_{X/Y}}$ and we call it a generalised Gysin map (or a virtual pull-back). Suppose furthermore that $Y$ is pure-dimensional. We define the virtual fundamental class of~$X$ with respect to the obstruction theory $E_{X/Y}$ to be the class
\[[X]_{E_{X/Y}}^{\vv}=f^!_{E_{X/Y}}[Y]=0^!_E[C_{X/Y}].\]
\end{construction}
\begin{Definition}We call $E_{X/Y}$ with property (2) an \emph{obstruction bundle} for $X\hookrightarrow Y$.
\end{Definition}
The following formula is useful for computing virtual classes.
\begin{Proposition}[{\cite[Example~4.1.8]{Fulton}}]\label{prop:zerocone}
Let $X\to Y$ be a closed embedding of schemes. If $C_{X/Y}$ is a cone of pure dimension $k$ and $E_{X/Y}$ a rank $r$ vector bundle containing it, we have
\begin{equation}\label{segrechern}0^!_E[C_{X/Y}]=\left\{c(E_{X/Y})\cap s(X,Y)\right\}_{k-r}=\sum_{i=0}^r{c_i(E_{X/Y})\cap s_{k-r+i}(C)},
\end{equation}
where $c(E_{X/Y})$ denotes the total Chern class of the vector bundle $E_{X/Y}$.
\end{Proposition}

\begin{Remark} For a more general definition of virtual classes see \cite[Definition 4.1]{siebert}. See \cite[Theorem~4.6]{siebert} for a generalisation of formula~\eqref{segrechern}.
\end{Remark}

\begin{Example}\label{excess}Let $X\hookrightarrow Y$ be a regular embedding with $Y$ pure-dimensional, then $C_{X/Y}=N_{X/Y}$ is a vector bundle. If $E_{X/Y}$ is an obstruction bundle, we may look at the following exact sequence of vector bundles
\[0\to N_{X/Y}\to E_{X/Y}\to E\to 0.\]
$E$ is called an \emph{excess bundle}. As the Segre class of a vector bundle is the inverse of its total Chern class (Example~\ref{ex:s=cinv}), from the above exact sequence we get that the virtual class is given by
\[[X]_{E_{X/Y}}^{\vv}=c_{\rm top}(E)\cap [X].\]
\end{Example}

\subsection{Obstruction bundles from global generators} \label{section:obs}
In order to define a virtual class on $X$, we need to embed $X$ in a pure-dimensional~$Y$ and find a~vector bundle~$E_{X/Y}$ which contains the normal cone~$C_{X/Y}$. Let us now construct such bundles.

Given an embedding $i\colon X\hookrightarrow Y$ of schemes and $f_1,\ldots, f_n$ generators of the ideal of~$X$ in $Y$, which are non-zero divisors, we have a surjective morphism
\[\bigoplus_{j=1}^rO_Y(-D_j)\xrightarrow{(f_1,\dots, f_r)}\mathcal{I}_{X/Y}\to 0, \]
where $D_j$ are Cartier divisors corresponding to $f_j$. Pulling this sequence back to $X$, we get an induced embedding of cones
\[C_{X/Y}\hookrightarrow\operatorname{Spec}\big({\mathcal S}{\rm ym}\big(\mathcal I_{X/Y}/\mathcal I_{X/Y}^2\big)\big)\hookrightarrow \operatorname{Spec}\left({\mathcal S}{\rm ym}\left(\bigoplus_1^ri^*\mathcal O_Y(-D_j)\right)\right).\]
We define an obstruction bundle for $i$
\[
E_{X/Y}:=\operatorname{Spec}\left({\mathcal S}{\rm ym}\left(\bigoplus_1^r\mathcal O_X(-i^*D_j)\right)\right).
\]
This is the vector bundle conventionally associated to $ \bigoplus_1^r\mathcal O_X(i^*D_j) $.

Note that $E_{X/Y}$ depends on the choice of $f_1,\ldots, f_n$ and thus also the virtual class in Construction~\ref{naive} depends on this choice.
\begin{Example}\label{doublepoint}
Let us compute the virtual class of $X=V\big(x^2,xy\big)\subset\mathbb{P}^2$, with respect to the vector bundle $E=\mathcal{O}_{\mathbb{P}^2}(2)\oplus \mathcal{O}_{\mathbb{P}^2}(2)$.

To see that this makes sense, note that $X=\{s=0\}$ for $s$ a section of $E$, which shows that $E^\vee\rvert_X\to I/I^2\to 0$. Let $0_E\colon X\to E$ denote the zero section of $E$. The virtual class of $X$ is defined by
\[[X]^{\vv}=0_E^!([C_{X/\mathbb{P}^2}]).
\]
By Example \ref{dorigin}, $[C_{X/\mathbb{P}^2}]=[C_0]+[L]$, where $C_0$ denotes the component supported on the origin and $L$ denotes the component of $C_{X/\mathbb{P}^2}$ supported on the $y$-axis. Since $0_E^!([C_{\mathbb{P}^2} X])=0_E^!([C_0])+0_E^!([L])$ we may deal with components of the cone separately.

Note that $C_0$ is the scheme given by the ring
\[
\mathbb C[x,y,A,B]/ \big(x^2,xy,y^2, yA-xB\big).
\]
Its reduced structure is the plane given by the ring $\mathbb C[A,B]$, and its multiplicity is two. This gives
\begin{equation}\label{first cone}
0_E^!([C_0])=0_E^!(\pi^*2[\underline 0])=2[\underline 0].
\end{equation}
We now compute $0_E^!([L])$. In the chart $z\neq 0$, where the singular point is,
\[L\subseteq C_{X/\mathbb{P}^2} \]
is given by $x=B=0$. This shows that $L$ is a line bundle locally defined by $B=0$. The normal bundle of the line $x=0$ is the line bundle $\mathcal{O}_{\PP^1}(1)$ and the line bundle $L$ is isomorphic to $\mathcal{O}_{\PP^1}(1)$ away from the origin. The embedding of $L$ in $\mathcal{O}_{\PP^1}(2)\oplus \mathcal{O}_{\PP^1}(2)$ is given by the injective morphism of sheaves
\[
0\to \mathcal{O}_{\PP^1}(1)\stackrel{(x,y)}{\longrightarrow} \mathcal{O}_{\PP^1}(2)\oplus \mathcal{O}_{\PP^1}(2).
\]
More precisely, $L$ is the closure of the image of $\mathcal{O}_{\PP^1}(1)$ inside $\mathcal{O}_{\PP^1}(2)\oplus \mathcal{O}_{\PP^1}(2)$.
Since on the line $x=0$, the above morphism is multiplication with $(0, y)$, we have that $L=\mathcal O_{\mathbb{P}^1}(2)$.
By Proposition \ref{prop:zerocone}, we obtain
\begin{equation}\label{second cone}
0_E^!([L])=\big\{c(E)\cap c(L)^{-1}\big\}_{0}=c_1(E/L)=2[P],
\end{equation}
where $[P]$ is the class of any point on $X^{\rm red}\simeq\mathbb{P}^1$. By (\ref{first cone}) and (\ref{second cone}), we have that
\[[X]^{\vv}=4[P].
\]
\end{Example}
\begin{Example}\label{exteriorpoint}
This is \cite[Example 4.2.2]{Fulton}. Consider a smooth surface $Y$, and let $A$, $B$ and~$D$ be effective Cartier divisors on $Y$, such that $A$ and $B$ meet transversally at a smooth point~$P$ away from $D$. Consider $A'=A+D$ and $B'=B+D$ and let $X$ be the scheme theoretic intersection of $A'$ and $B'$. Let us compute the Segre class and the virtual class of~$X$ induced by $E_{X/Y}=\mathcal O_Y(A')\oplus\mathcal O_Y(B')$. Note that $A'$ and $B'$ correspond to generators for $\mathcal{I}_{X/Y},$ and $E_{X/Y}$ is an obstruction bundle for $X\to Y$.

We compute the Segre class after blowing up, as in Proposition~\ref{prop2}. Let $\tilde Y$ be the blow-up of~$Y$ at $P$, let $\tilde D$ be the strict transform of $D$ and let $E$ be the exceptional divisor. $\tilde Y$ is isomorphic to the blow up of $Y$ along $X$ and under this isomorphism $\tilde X$ corresponds to $\tilde D+E$. Recall that $s(X,Y)=\pi_*\big(\big[\tilde X\big]-\tilde X\cdot\big[\tilde X\big]\big)$, hence we need to compute $\big(\tilde D+E\big)^2$. Using the formulas in Lemma~\ref{blup_formulas}, we obtain
\begin{equation*}
 \tilde D^2=D^2 \qquad\text{and} \qquad \tilde D\cdot E=0 \qquad\text{and}\qquad E^2=-[P_E],
\end{equation*}
where $P_E$ is the class of a point on $E$. This is because the blow-up does not change a neighbourhood of $D$ and $\tilde D$ and~$E$ do not intersect. Then,
\[
 s(X,Y)=\pi_*\big(\big[\tilde D+E\big]-\big(\tilde D+E\big)\cdot\big[\tilde D+E\big]\big)=[D]-D\cdot[D]+[P].
\]
The virtual class on $X$ is given by
\begin{align*}
 [X]^\text{vir}&=\{c(E_{X/Y})\cap s(X,Y)\}_0 =c_1(E_{X/Y})\cap[D]-D\cdot[D]+[P]
 \\&=(D+A+B)\cdot[D]+[P].
\end{align*}
\end{Example}
\begin{Example}Let us now look at a smooth surface $Y$, and $A$, $B$, $D$ effective Cartier divisors on $Y$ as before, but now let $A$ and $B$ meet transversally at a smooth point $P\in D$. For example, consider $Y=\PP^2$ with homogeneous coordinates $[x:y:z]$, $A=\{y=0\}$, $B=\{x-y=0\}$, $D=\{x=0\}$. We then have that the ideal of $A'\cap B'$ is $(xy, x(x-y))$ or equivalently $\big(xy, x^2\big)$. This shows that $A'\cap B'$ is isomorphic to $X$ in Example~\ref{doublepoint}. Although the set theoretic intersection of $A'$ and $B'$ is $D$, the scheme theoretic intersection is different and the virtual class is different from $D$. This is compatible with the fact that we can deform $A$ and $B$ in Example~\ref{exteriorpoint} such that in the limit $P\in D$ and virtual classes are deformation invariant. We can easily check deformation invariance in this case. If we take $A$, $B$, $D$ in Example~\ref{exteriorpoint} to be lines, we get $[X]^{\vv}= 4[P]$. This agrees with Example~\ref{doublepoint}.
\end{Example}

\section{Gluing obstruction bundles}\label{gluingobs}
In this section we motivate and define obstruction theories. We give the general definition of virtual classes as it appears in \cite{BerFan}. For many aspects in \cite{BerFan} which we do not discuss here, we refer to the original paper by Behrend and Fantechi. We refer to~\cite{Kstacks} and~\cite{vir-pull} for details regarding the stacky version of virtual pull-backs.

For what we have seen so far, it is often easy to define a virtual class on a scheme. What is difficult is to define a virtual class of a given \emph{expected dimension} (see Section~\ref{obsdefth}). Recall that the dimension of the virtual cycle $[X]_{E_{X/Y}}^{\vv}$ depends on the rank of the obstruction bundle~$E_{X/Y}$. Obstruction bundles induced by generators of an ideal as in Section~\ref{section:obs} may often lead to trivial virtual classes. It is then natural to look for obstruction bundles of lower rank. Let us see what can go wrong in an example.
\begin{Example}
Consider the embedding of the twisted cubic
\[
 \begin{tikzcd}[row sep=tiny]
 \mathbb P^1\ar[r] & \mathbb P^3, \\
 \left[s:t\right]\ar[r, mapsto] & \big[s^3:s^2t:st^2:t^3\big]
 \end{tikzcd}
\]
with ideal $I_{X/\PP^3}$ generated by $\big(xz-y^2,yw-z^2,xw-yz\big)$. Let $A=xz-y^2$, $B=yw-z^2$, $C=xw-yz$. According to the construction in Section~\ref{section:obs}, the obstruction bundle restricted to~$X$ is $E_{X/\PP^3}=\oplus_1^3 \mathcal O_{\PP^1}(6)$. On the other hand, $X$ is smooth, so we can take the normal bundle as obstruction bundle. In the following, we construct the normal bundle from equations and gluing data. We see that we have two relations among the generators of $I_{X/\PP^3}$. Let us pick one of them, for example $wA+By-zC=0$. This relation induces a morphism
\[ \oplus_1^3 \mathcal O_{\PP^1}(6)\to \mathcal O_{\PP^1}(9),
\]
given by diagonal multiplication with $(w,y,-z)$. The kernel of this morphism is $N_{X/\PP^3}$. The cokernel is $\mathbb{C}_P$, where $P=[1:0:0:0]$. This shows that we have an exact sequence
\[
0\to N_{X/\PP^3}\to \oplus_1^3 \mathcal O_{\PP^1}(6)\to \mathcal O_{\PP^1}(8)\to 0.
\]
One can check that $N_{X/\PP^3}\simeq \mathcal O_{\PP^1}(5)\oplus\mathcal O_{\PP^1}(5)$, but we will not need this here. By the excess intersection formula in Example \ref{excess}, the virtual class with respect to $E_{X/\PP^3}$ is
 \[
 [X]^{\vv}_{E_{X/\PP^3}}=8[pt].
 \]
 Taking the obstruction bundle of $X$ in $\PP^3$ equal to $N_{X/\PP^3}$, we get $[X]^{\vv}=[X]$.
\end{Example}
 \begin{Example}\label{gluecone} Let us look at one scheme which is not smooth. In this case an obstruction bundle is not obvious. Consider $\PP^{5}\times \PP^4$ with coordinates $\left([y_0,\ldots, y_4,\xi],[x_0,\ldots,x_4]\right)$ and $X\subset \PP^5\times\PP^4$, with equations
 \begin{gather*}y_0\xi,\ldots, y_4\xi,\\
 x^s_0y_0+\dots +x^s_4y_4,\\
 Q_1,\ldots,Q_l
 \end{gather*}
 with $Q_j$ homogeneous polynomials in $x_0,\ldots, x_{4}$. Since we do not impose any conditions on~$Q_i$, $X$~can be very singular and may have many components, so the cone $C_{X/\PP^5\times\PP^4}$ may not be a~bundle as we had in the previous example. As before, we have an obstruction bundle $E$ given by global generators $A_i=\xi y_i$, $B=x^s_0y_0+\dots+ x^s_4y_4$, $C_j=Q_j$, but we can find a smaller one in the following way. Let $L$ be the line bundle $\OO_{\PP^5\times\PP^4}(2,s).$ Then we have a morphism
 \begin{align*}
 E&\to L, \\
 (A,B,C)&\mapsto x^s_0A_0+\dots+x^s_4A_4-\xi B.
 \end{align*}

 Since not all $x_i$ are zero, this morphism is surjective. Let $E'$ denote the kernel of $ E\to L$. Since the equation $x^s_0A_0+\dots+x^s_4A_4-\xi B$ vanishes on $C_{X/\PP^5\times\PP^4}$, we have that $E'$ is an obstruction bundle for $X\to \PP^5\times\PP^4$. $E'$ is sometimes called a reduced obstruction theory.
 \end{Example}

 Let us now give a new definition of obstruction bundles motivated by the previous examples. The following definition is less coarse than the one in Section~\ref{section:obs}, but still rather na\"ive.
 \begin{Definition}\label{obs: less coarse}Let $X$ be a DM stack (or a scheme). An obstruction theory $E_{X/V}$ on $X$ is given by the following data:
\begin{enumerate}\itemsep=0pt
\item[1)] a global embedding $X\hookrightarrow V$, with $V$ a smooth stack,
\item[2)] a covering of $V$ by open sets $V_i$, and on each $U_i=X\cap V_i$ a set of generators for the restriction to $U_i$ of the ideal $I_{U_i/V_i}$; these induce obstruction bundles~$E_i$,
\item[3)] gluing data for the sets of generators, i.e., isomorphisms
\[\phi_{ij}\colon \ E_i|_{U_{ij}}\to E_j|_{U_{ij}}
\]
such that $\phi_{ij}$ satisfy cocycle and such that the following diagram commutes
\[\xymatrix{C_{U_i/V_i}|_{U_{ij}}\ar[r]\ar[d]_{\bar{\phi}_{ij}}&E_i|_{U_{ij}}\ar[d]^{\phi_{ij}}\\
C_{U_j/V_j}|_{U_{ij}}\ar[r]&E_j|_{U_{ij}}.}
\]
Here $\bar{\phi}_{ij}$ is the canonical isomorphism, and the horizontal maps are the embeddings of the cones in the corresponding obstruction theories.
\end{enumerate}

 \end{Definition}
\begin{Remark} Given $X$ in $V$ as above and equations of $X$ in $V$, everything in the definition above is explicit. We have seen that cones can be described explicitly inside the corresponding obstruction bundles and gluing data can be given in terms of relations between generators of~$I_{U_i/V_i}$ and~$I_{U_j/V_j}$.
\end{Remark}
In conclusion, if we want to get reasonable virtual cycles, we need a lot of information: local equations of our space in a smooth ambient space, and gluing data for local obstruction bundles. In this way computations get quickly very messy. The information we want can be encoded more efficiently in a stacky gadget that we describe below. For this, we remark that given the embedding $X\hookrightarrow Y$, the na\"ive obstruction bundle $E_{X/Y}=\bigoplus \mathcal O_X(D_i)$ comes with a~piece of extra data we did not use so far. We have a morphism
\begin{equation}\label{action}T_Y\rvert_X\xrightarrow{(\delta f_i)} E_{X/Y},
\end{equation}
defining an action of the tangent space of the ambient restricted to $X$. This motivates the following.

Let $X$ be a scheme (DM stack). In the following $\mathbb L^\bullet_{X}\in D^{\leq 0}(X)$ denotes the cotangent complex of $X$, and $\mathbb L^{\geq -1}_{X/Y}$ the $[-1,0]$ truncation of $\mathbb L^\bullet_{X}$. See Remark~\ref{cotcpx} for examples of truncated cotangent complexes.

\begin{Definition}[{\cite[Definition 4.4]{BerFan}}]\label{bfobs}
A \textit{perfect obstruction theory} for $X$ is a complex $E^{\bullet}_{X}\in D^{b}(X)$ which is locally quasi-isomorphic to a complex of locally free sheaves supported in $[-1, 0]$, with a morphism $E^{\bullet}_{X}\xrightarrow{\phi}\mathbb L^{\geq -1}_{X} $ such that
 \begin{enumerate}\itemsep=0pt
 \item[1)] $h^0(\phi)\colon h^0(E^{\bullet}_{X})\to h^0\big(\mathbb L^{\geq -1}_{X}\big)$ is an isomorphism,
 \item[2)] $h^{-1}(\phi)\colon h^{-1}(E^{\bullet}_{X})\to h^{-1}\big(\mathbb L^{\geq -1}_{X}\big)$ is surjective.
\end{enumerate}
 \end{Definition}
In particular, we get a perfect obstruction theory if locally on $X$ we may find embeddings $U\subseteq V$, with $V$ smooth, together with obstruction bundles~$E_{U/V}$, such that the local tangent-obstruction complexes $[T_V|_U\to E_{U/V}]$ glue to $E^{\bullet}_{X}$, which is a global object on $X$. Such a~complex gives rise to an algebraic stack denoted $h^1/h^0(E^{\bullet}_{X})$ such that \begin{gather*}h^1/h^0(E^{\bullet}_{X})|_U\simeq E_{U/V}/T_V|_U.\end{gather*} See \cite[Section~2]{BerFan} for the definition. We use the notation $\mathfrak{E}_X=h^1/h^0(E^{\bullet}_{X})$. We call $\E_X$ a~\emph{vector bundle stack}. The rank of a vector bundle stack $h^1/h^0(E^{\bullet}_{X})$ as above, is defined to be $\rk\big(h^1(E^{\bullet}_{X})\big)-\rk\big(h^0(E^{\bullet}_{X})\big)$. Notice that this number may be negative.

We have an analogous construction for cones. Given local charts $U_i$ of our space $X$ which embed into possibly unrelated smooth ambient spaces $V_i$, the action in~(\ref{action}) preserves the cone~$C_{U/V}$ and we have a stack $\mathfrak C_X$ such that
\[\mathfrak C_X|_U\simeq C_{U/V}/T_V|_U.
\]
The stack $\mathfrak C_X$ is called \emph{the intrinsic normal cone} of~$X$. It is a zero dimensional stack, independent of embeddings in smooth stacks. See~\cite{BerFan} for the construction.

Conditions (1) and (2) in Definition~\ref{bfobs} are equivalent to the fact that the vector bundle stack $\mathfrak{E}_X$ associated to $E^{\bullet}_{X}$ contains the intrinsic normal cone $\mathfrak C_X$ \cite[Theorem~4.5]{BerFan}.

With these in place, we can now define a more general virtual class.
 \begin{Proposition}[\cite{Kstacks}] If $\E$ is a vector bundle stack on a scheme $($or DM stack$)$~$X$, the flat pull-back is an isomorphism.
\end{Proposition}
We denote its inverse by
\[0^!_{\E}\colon \ A_*(\E)\xrightarrow{\cong} A_{*-\rk(\mathfrak{E})}(X).\]

\begin{construction}[\cite{BerFan}] \label{bfgen} Given a perfect obstruction theory $E^\bullet_X$ on a DM stack $X$, we may define its virtual class by means of
\[[X]^{\vv}_{\E_{X}}=0^!_{\E_{X}}[\mathfrak C_X].\]
\end{construction}
\begin{Remark} Unlike Construction~\ref{naive}, the construction in~\cite{BerFan} only requires local embeddings in a smooth stack. This is a major advantage over Construction~\ref{naive}, as global embeddings in smooth stacks are sometimes unavailable. However, Construction~\ref{bfgen} is in general impractical for computational purposes. For this reason, in this note we will mainly work with the more na\"ive Construction~\ref{naive}.
\end{Remark}
\begin{Remark}
 Construction \ref{bfgen} may be generalised to any morphism of DM type $f\colon X\to Y$ between algebraic stacks. The relative intrinsic normal cone may locally be constructed as follows: choose a~chart~$U$ for $X$, embedded in a space $V$ smooth over $Y$ (e.g., embed $U$ in a~smooth~$M$ and take $V=Y\times M$), then $\mathfrak C_{X/Y}|_U=[C_{U/V}/T_{V/Y}|_U]$. This leads to the definition of virtual pull-back~\cite{vir-pull}
 \[f^!_{\E_{X/Y}}\colon \ A_*(Y)\to A_{*-\rk(\E_{X/Y})}(X).\]
 It naturally endows $X$ with a virtual class whenever $Y$ is pure-dimensional.
\end{Remark}

\begin{Remark} We are usually interested in finding obstruction theories for moduli functors. In many cases obstruction theories are available from deformation theory. More precisely, consider the DM stack~$X^F$ representing a certain moduli functor~$F$. In this case, the study of the deformation theory of $F$ always gives a local tangent-obstruction theory, i.e., the Kuranishi map:
\[T^1_F\to T^2_F.\]
A perfect obstruction theory is usually given considering a global version of this picture, see the examples in~\cite{Ber}. Similar constructions can be done for other moduli spaces (quasimap spaces~\cite{CF-K-M}, certain moduli spaces of sheaves~\cite{thomas} or stable pairs on Calabi--Yau 3-folds~\cite{pan-thomas}). See~\cite{rt} for an overview.
\end{Remark}

\begin{Remark} Above we have constructed $\E_X$ from a two-term complex $E^{\bullet}_X\in D^{[-1,0]}(X)$. There is a more general construction due to Chang and Li (see~\cite{cl:semiperfect}). Their ingredients are a~semi-perfect obstruction theory and a Gysin map from the obstruction sheaf. A semi-perfect obstruction theory consists of local obstruction theories and gluing data for the obstruction sheaf.
\end{Remark}

\begin{Remark}\label{cotcpx} We conclude this section by describing the cotangent complex in some simple cases.
\begin{enumerate}\itemsep=0pt
\item If $X\hookrightarrow Y$ is an embedding, then
\[\mathbb L^{\geq -1}_{X/Y}=[\mathcal I/\mathcal I^2\to 0].\]
In this case, any vector bundle $E$ which surjects onto $\mathcal I/\mathcal I^2$ gives a perfect obstruction theory. Many of the examples in this note are of this type.
\item Let $Y$ be a point and $X$ a smooth DM stack, then $\mathbb L^\bullet_{X/pt}\cong \Omega_X$ and given any vector bundle~$E$, the complex $\big[E\xrightarrow{0} \Omega_X\big]$
is a perfect obstruction theory. Notice that the resulting virtual class is $[X]^{\vv}=c_{\text{top}}(E)\cap [X]$.
\item Let $g$ be an embedding of $X\hookrightarrow Y\times M$ relative to $Y$, with $M$ smooth, then the truncated relative cotangent complex has the following form
\[\mathbb L^{\geq -1}_{X/Y}=\big[\mathcal I/\mathcal I^2\to\Omega_{M\times Y/Y}\big],\]
where $\mathcal I$ is the ideal of the embedding. If $X$ is l.c.i.\ in $Y\times M$ then the above is a two-term complex of vector bundles, and it represents the full cotangent complex $\mathbb L^\bullet_{X/Y}$.
\end{enumerate}
\end{Remark}

\section{Deformations of obstructions and cones}\label{deformations}
In this section we discuss functoriality of virtual pull-backs. We give some elementary examples. The main references for the general theory are~\cite{KKP} and~\cite{vir-pull}.

\begin{construction}[deformation of vector bundles]\label{defvb}
Let
\[
0\rightarrow E\xrightarrow{\phi} F\rightarrow G\rightarrow 0
\]
be a sequence of vector bundles on a scheme/stack~$X$. Consider the family of morphisms
\[
0\rightarrow E\xrightarrow{(t\cdot id_E,\phi)} E\oplus F
\]
on $X\times\mathbb A^1$. The Cokernel $\operatorname{coker}((t\cdot {\rm id}_E,\phi)$ is a family of vector bundles
\begin{align*}
\xymatrix{\mathcal{E}\ar[d]^{\pi}\\
\mathbb{A}^1}\;\;\;\text{s.t.}\;\;\mathcal{E}\rvert_{\mathbb{A}^1\setminus\{0\}}\cong F\times \mathbb{A}^1\setminus\{0\}\quad
\text{and}\quad \pi^{-1}(0)=E\oplus G.
\end{align*}
Given a cone $C\subseteq F$, considering the closure of $C\times \mathbb{A}^1\setminus \{0 \}$ in $\E$, we get a central fibre $C_0$ which we call the limit of $C$ in $E\oplus G$. By construction, the limit $C_0$ is rationally equivalent to~$C$.
\end{construction}

\begin{Proposition}\label{easycompat}
Given embeddings $X\hookrightarrow Y\hookrightarrow Z$ and an exact sequence of obstruction bundles
\[0\to E_{X/Y}\to E_{X/Z}\to E_{Y/Z}\rvert_X\to 0\]
compatible with the maps from the normal sheaves,
the virtual classes on $X$ given by
\[C_{X/Z}\hookrightarrow E_{X/Z}\qquad\text{and}\qquad C_{X/C_{Y/Z}}\hookrightarrow E_{X/Y}\oplus E_{Y/Z}\rvert_X\]
coincide.
\end{Proposition}
\begin{proof}[Proof (sketch)]
Consider first the deformation to the normal cone for the map $Y\hookrightarrow Z$. This is a flat one-parameter family $M^0_{Y/Z}\to\PP_{ t}^1$, such that the fibers over $t\neq 0$ are isomorphic to~$Z$, and the fiber over~$0$ is isomorphic to $C_{Y/Z}$ (see Section~\ref{sec:defnormal}).
Construct now the double deformation space
\[
M_{s,t}=M^0_{{X\times\PP^1_t}/M^0_{Y/Z}}\to \PP^1_t\times \PP^1_s.
\]
This is a flat family over $\PP^1_s$, with general fibre over $s\neq 0$, the deformation space $M^0_{Y/Z}$ and special fibre over $s=0$, the cone $C_{X\times \PP^1_t/M^0_{Y/Z}}$. One can check that the fibres of $M_{s,t}$ over $ \PP^1_t\times \PP^1_s$ are as follows:
\begin{enumerate}\itemsep=0pt
\item[1)] the fibre over $s,t\neq 0$ is $Z$;
\item[2)] the family over $s\neq 0$ fixed, when $t$ varies is $M_{s_0,t}=M^0_{Y/Z}$; this family is flat over $\PP^1_t$;
\item[3)] the family over $s=0$, when $t$ varies is the one parameter family $M_{t}=C_{X\times \PP^1_t/M^0_{Y/Z}}$; this family over $\PP^1_t$ is not flat;
\item[4)] the family over $t\neq 0$, when $s$ varies is $M_{s,t_0}=M^0_{X/Z}$; this family is flat over $\PP^1_s$;
\item[5)] the family over $t=0$, when $s$ varies is the one parameter family $M_{s}$; $M_s$ contains the flat family $M^0_{X/C_{Y/Z}}\to \PP^1_s$ and it agrees with it generically; this family is not flat.
\end{enumerate}
Let $t_0\colon 0\to \PP^1_t$ be the embedding of the point $0$ and $s_0\colon 0\to \PP^1_s$ be the embedding of the point~$0$.

Let us now look at $M_t=C_{X\times \PP^1_t/M^0_{Y/Z}}\to \PP^1_t$. We have that
\begin{gather*}t_0^![M_t]=t_1^![M_t]
\end{gather*}
in $A_*(M_t)$. Since the general fibre of $M_t\to\PP^1_t$ is $C_{X/Z}$, we have that
\begin{equation*}
t^!_1[M_t]=[C_{X/Z}].
\end{equation*}

In the following we show that $t_0^![M_t]=[C_{X/C_{Y/Z}}]$.\footnote{The fibre over $t=0$ contains $C_{X/C_{Y/Z}}$, but it may be larger in general, so we need a more elaborate argument.}

By commutativity of intersection products, we have that
\begin{equation}\label{s comm t} s_0^!t^!_0[M_{s,t}] =t_0^!s^!_0[M_{s,t}]
\end{equation}
in $A_*(M_{0,0})$.

One can see from the construction of $M_{s,t}$, and point~(2) in the description of the fibres of~$M_{s,t}$ above, that
\begin{equation}\label{restriction1}t^!_0[M_{s,t}]=[M^{\prime}_s] \qquad \text{and} \qquad s^!_0[M_{s,t}]=[M_t].
\end{equation}
Here $M^{\prime}_s$ is obtained by discarding components of $M_{s,t}$ which do not dominate $\PP^1_t$, followed by restriction to $t=0$. Since the map $M_{s,t}\to\PP^1_t$ is generically trivial, such components are embedded in the fibre over $t=0$. Suppose now that we have components of $M_{s,t}$ in the fibre $t=0$, but not over $(s,t)=(0,0)$. Then, we have that a generic fibre $M_{s,t_0}$, with $t_0\neq 0$, contains the restriction of such a component. This contradicts point (4) in the above description of fibres of $M_{s,t}$. Since $s_0^!$ discards components of $M_{s,t}$ which map to $s=0$, the above discussion shows that
\begin{equation}\label{restriction2}s^!_0[M_s]=s^!_0[M^{\prime}_s].
\end{equation}
Since $M_s$ contains $M^0_{X/C_{Y/Z}}$, which is flat over $\PP^1_s$, we have that
\begin{equation}\label{sdirflat}s^!_0[M_s] =s^!_0[M^0_{X/C_{Y/Z}}]=[C_{X/C_{Y_Z}}]
\end{equation}
in $A_*(M_{0,0})$.
\noindent By (\ref{s comm t}), (\ref{restriction1}), (\ref{restriction2}) and (\ref{sdirflat}), we have that $t_0^![M_t]=[C_{X/C_{Y/Z}}]$ in $A_*(M_{0,0})$. This shows that we have a rational equivalence
\begin{equation}\label{equivcones}[C_{X/Z}]\simeq [C_{X/C_{Y/Z}}]
\end{equation}
in $A_*(C_{X\times\PP^1_t/M^0_{Y/Z}})$.
The hypothesis implies that we have an embedding
\[
M_t\hookrightarrow \mathcal{E},
\]
where $\mathcal{E}$ is the deformation of obstruction bundles $E_{X/Z}\rightsquigarrow E_{X/Y}\oplus E_{Y/Z}\rvert _X$ (see \cite{KKP}).
Pushing this equivalence forward to $\mathcal{E}$, we get that \eqref{equivcones} holds in $A_*(\mathcal{E})$. The claim on virtual classes follows from commutativity of $0^!_{\mathcal E}$ and restriction to $t\in\PP^1$.
\end{proof}

\begin{Example}\label{mainexample}
Let us compute the virtual class of $X=V(xz,yz)\subseteq\PP^3$. The embedding of~$X$ in~$\PP^3$ factors through $Y=V(xz)$. This gives an exact sequence of compatible obstruction theories:
\[\xymatrix{0\ar[r] &\OO_X(2)\ar[r]\ar@{=}[d] &\OO_X(2)\oplus\OO_X(2)\ar[r]\ar@{=}[d] &\OO_X(2)\ar[r]\ar@{=}[d] &0.\\
& E_{X/Y} &E_{X/\PP^3} &E_{Y/\PP^3}\rvert_X}\]
By Proposition~\ref{easycompat}, the virtual class is given by the embedding
\[C_{X/C_{Y/\PP^3}}\rvert_X\hookrightarrow E_{X/Y}\oplus E_{Y/\PP^3}.\]
More precisely, we have
\[[X]^{\vv}=0^!_{E_{X/Y}\oplus E_{Y/\PP^3}}[C_{X/C_{Y/\PP^3}}].
\]
Note that $C_{Y/\PP^3}\simeq\OO_Y(2)$, so
\[
C_{X/C_{Y/\PP^3}}\cong C_{X/Y}\times_X C_{Y/\PP^3}\rvert_X \cong C_{X/Y}\oplus \OO_X(2).
\]
Moreover, as $Y=Y_1\cup Y_2=V(z)\cup V(x)$
\[C_{X/Y}\cong C_{X_1/Y_1} \cup C_{X_2/Y_2}=C_{H/H}\oplus C_{L_1\cup L_2/Y_1},\]
where \begin{gather*}X_1=X\times_Y Y_1=H=V(z)\end{gather*} and
\begin{gather*}X_2=X\times_Y Y_2=L_1\cup L_2=V(x,y)\cup V(x,z).\end{gather*}
We now compute the virtual class of $X$ with respect to the obstruction theory $E_{X/\PP^3}=\OO_X(2)\oplus \OO_X(2)$. The two components of $C_{X/C_{Y/\PP^3}}$ give the following contributions.
\begin{enumerate}\itemsep=0pt
\item From the first component of $C_{Y/\PP^3}$ we get a contribution
\[C_{X_1/Y_1}\times_{X_1}\OO_{Y}(2)\rvert_{X_1}=\OO_H(2)\hookrightarrow \OO_H(2)\oplus \OO_H(2)=E_{X/\PP^3}\rvert_{X_1} .\]
 This gives $2[L]$, where $L$ is the class of a line on the plane $H=V(z)$;
\item From the second component of $C_{Y/\PP^3}$ we get the contribution
\begin{align*}C_{L_1\cup L_2/Y_1}\times_{X_2}\OO_Y(2)\rvert_{X_2}\cong \OO_{L_1\cup L_2}(2)^{\oplus 2}\hookrightarrow \OO_{L_1\cup L_2}(2)^{\oplus 2} =E_{X/\PP^3}\rvert_{X_2}
\end{align*}
 We thus obtain a contribution $[L_1]+[L_2]$.
\end{enumerate}
In conclusion $[X]^{\vv}=2[L]+[L_1]+[L_2]$. We can re-write this as
\[[X]^{\vv}=[L_1]+3[L].\]
\end{Example}

\begin{Example}\label{limit main example}Let us now compute the limit $\mathcal{C}_0$ defined in Construction \ref{defvb}. By Example \ref{main example cone} we have an embedding
\[
C_{X/\PP^3}\hookrightarrow\OO_X(2)\oplus\OO_X(2).
\]
We apply Construction \ref{defvb} to the embedding
\[\iota_1\colon \ \OO_X(2)\hookrightarrow \OO_X(2)\oplus\OO_X(2).\]
Looking at coordinate rings we get that the induced $\iota_1$ is
\begin{align*}
R[A,B]&\twoheadrightarrow R[C],\\
A&\rightarrow C,\\
B&\rightarrow 0.
\end{align*}
We now deform $\iota_1$ to a family of embeddings over $\mathbb A^1_t$. Consider the morphism in Construction~\ref{defvb}
\[\OO_X(2)\xrightarrow{(\iota_1,t\cdot {\rm id})} \OO_X(2)\oplus\OO_X(2)\oplus\OO_X(2)\rightarrow\OO_X(2)\oplus\OO_X(2).\]
For $t\neq 0$, this corresponds to the morphisms of coordinate rings
\begin{align*}
\xymatrix{0\ar[r] & R[A,B]\ar[rr]^{\left(\substack{A\to A-\frac{1}{t}A'\\ B\to B}\right)} & & R[A',A,B]\ar[rr]^{\left(\substack{A'\to tC\\ A\to C\\ B\to 0}\right)} & & R[C]\ar[r] & 0.}
\end{align*}
Let $\mathcal C\hookrightarrow \mathcal{E}$ be the flat family in Construction~\ref{defvb}. By Example \ref{main example cone} we have that the normal cone~$C_{X/\PP^3}$ is locally isomorphic to $\operatorname{Spec}R[A,B]/(yA-xB)$.
This shows that $\mathcal C\hookrightarrow \mathcal{E}$ is locally given by $R[A',A,B]\twoheadrightarrow R[A',A,B]/(y(A-\frac{1}{t}A')-xB)$.
Taking the flat limit for $t$ goes to $0$, we get
\[\mathcal C_0\cong R[A',B]/yA'.\]
We see that in this case $\mathcal C_0= C_{X/C_{Y/Z}}$.
\end{Example}

\begin{Example}\label{nonlimit}In general, given embeddings $X\hookrightarrow Y\hookrightarrow Z$ and an embedding of $E_{X/Y}$ into an obstruction bundle $E_{X/Z}$ for $X\hookrightarrow Z$, the deformation $\mathcal C$ of the normal cone inside the family of cokernels does not necessarily give $\mathcal C_0\cong C_{X/C_{Y/Z}}$ (see Proposition~\ref{easycompat} and Example~\ref{limit main example}).
Consider
 \[
 X=V(x,y)\hookrightarrow Y=V(xy)\hookrightarrow Z=\PP^3_{[x:y:z:w]}.
 \]
We have obstruction theories $E_{X/Y}=\OO_X(1)\oplus\OO_X(1)$ and $E_{Y/Z}=\OO_Y(2)$. Set
\[
E_{X/\PP^3}=\OO_X(1)\oplus\OO_X(1)\oplus\OO_X(2).
\]

We now look at the deformation induced by the exact sequence of vector bundles
 \[
 0\to \OO_X(1)\oplus\OO_X(1)\to\OO_X(1)\oplus\OO_X(1)\oplus\OO_X(2)\to\OO_X(2)\to 0.
 \]
Let us compute the double deformation space. By restricting to the chart $\{w\neq 0\}$, we see that $\operatorname{Bl}_Y Z$ inside $\operatorname{Bl}_{Y\times\{0\}} Z\times\mathbb A^1_t\subseteq \mathbb A^4_{x,y,z,t}\times\PP^1_{[\alpha_0:\alpha_1]}$ is given by $\alpha_1=t=0$. This shows that its complement $M^0_{Y/Z}$ is the subscheme $V(xy-t\alpha_0)\subseteq \mathbb A^5_{x,y,z,t,\alpha_0}$, and $X\times\mathbb A^1_t$ is $V(x,y,\alpha_0)$. Similarly, we find that the double deformation space $M^0_{X\times\mathbb A^1_t/M^0_{Y/Z}}$ inside $\mathbb A^6_{z,t,s,\beta_0,\beta_1\beta_2}$ has equation $\{s\beta_0\beta_1-t\beta_2=0\}$. Here $x=s\beta_0$, $y=s\beta_1$, $\alpha_0=s\beta_2$. Restricting to $s=0$, we see that the flat limit of $C_{X/Z}\times\{t\neq 0\}$ is $C_{X/Z}$. On the other hand, if we restrict to $t=0$, we find a cone with two components. More precisely, it agrees with the following
\begin{align*}
C_{X/C_{Y/Z}}&=C_{X/Y}\times_X C_{Y/Z}\rvert_X =C_{X/H_1}\times_{X} C_{Y/Z}\rvert_{X}\cup C_{X/H_2}\times_{X} C_{Y/Z}\rvert_{X}\\
&=\OO_X(1)\oplus\OO_X(2)\cup \OO_X(1)\oplus\OO_X(2).
\end{align*}

Notice that $C_{X/Z}$ and $C_{X/C_{Y/Z}}$ have different Segre classes, in particular they cannot lie in a~flat family of cones over $X$ (see Proposition \ref{Segre constant}). Nevertheless, they are rationally equivalent in the obstruction bundle, hence the virtual class computation produces the same result
\[0^!_{E_{X/Z}}([C_{X/Z}])=c_1(\OO_X(2))=2 c_1(\OO_X(1))=0^!_{E_{X/Y}\oplus E_{Y/Z}\rvert_X}([C_{X/C_{Y/Z}}]).
\]
\end{Example}

\begin{Remark}Let us look at $\mathcal C_{X\times\mathbb A^1/M_{Y/Z}^\circ}$. In a first step, we consider $\operatorname{Bl}_{Y\times\{0\}}Z\times\mathbb A^1_t$. Setting $f=xy$, this is \begin{gather*}\Proj \mathbb C[t][x,y,z][\alpha,\beta]/(t\alpha-f\beta).\end{gather*} We have to remove $\operatorname{Bl}_{Y}Z$ from the central fiber, that is the locus where $\beta=0$. One sees that $X\times \mathbb A^1$ is given by equations $\{x=y=\alpha=0\}$. With this, we get that the normal cone $\mathcal C_{X\times\mathbb A^1/M_{Y/Z}^\circ}$ is given by the ring $\mathbb C[t][w_0,w_1,w_2]/(tw_2)$. This has two components, one of them being supported on $\{t=0\}$. This shows that $\mathcal C_{X\times\mathbb A^1/M_{Y/Z}^\circ}$ is not flat over $\mathbb A^1_t$.
\end{Remark}

\section{Properties of virtual classes and computational methods}\label{properties}

Virtual classes have good functorial properties which we now list. We then discuss push-forwards of virtual classes.

\subsection{Pull-backs of virtual classes}
\begin{Theorem}[see \cite{vir-pull} for the general statement and {\cite[Chapter 6]{Fulton}} for regular embeddings]\label{vir pull-back} Let $f\colon X\to Y$ be a morphism of DM type with a perfect obstruction theory $E^\bullet_{X/Y}$. Then $f^!_{\mathfrak E_{X/Y}}$ is a bivariant class $($see Chapter~$17$ in~{\rm \cite{Fulton})}.

In particular, the following properties hold.
\begin{enumerate}\itemsep=0pt
\item[$(a)$] Given a Cartesian diagram
\[
\xymatrix{X'\ar[r]^g\ar[d]_p&Y'\ar[d]^q\\
X\ar[r]^f&Y}\]
with $q$ proper such that $E^\bullet_{X/Y}$ is a perfect obstruction theory for $f$, then $p^*E^\bullet_{X/Y}$ is an obstruction theory for $g$, and for any $\alpha\in A_*(Y')$ we have that
\[p_*g^!_{p^*\mathfrak E}\alpha=f^!_{\mathfrak E}q_*\alpha
\]
in $A_*(X)$.

\item[$(b)$] Given a Cartesian diagram
\[
\xymatrix{X'\ar[r]\ar[d]&Y'\ar[d]^g\\
X\ar[r]^f&Y}\]
such that $E^\bullet_{X/Y}$ is a perfect obstruction theory for $f$, and $E^\bullet_{Y'/Y}$ is a perfect obstruction theory for $g$, then for any $\alpha\in A_*(Y)$, we have that
\[g^!f^!\alpha=f^!g^!\alpha
\]
in $A_*(X')$.
\end{enumerate}
\end{Theorem}
\begin{Remark}Let us look again at the line with an embedded point in Example \ref{doublepoint}.
By Theo\-rem~\ref{vir pull-back} and Example~\ref{excess}, applied to the Cartesian diagram
\[\xymatrix{X\ar[r]^i\ar[d]&\PP^2\ar[d]^{x^2,xy}\\
\PP^2\ar[r]& E
}
\]
we have that
\begin{gather*}i_*[X]^\text{vir}= c_\text{top} (E)\cap \big[\mathbb P^2\big] =4[pt],
\end{gather*}
where $pt$ is a point in $\PP^1$. By looking at the components of the normal cone, as we did in Example~\ref{doublepoint}, we saw that two of those points are supported on the origin. This is however, irrelevant in this example; since
\[A_0(X)=A_0\big(X^{{\rm red}}\big)=\mathbb Z,
\]
we have that all points of $X$ are rationally equivalent.

Example \ref{mainexample} -- a line and a plane -- is more meaningful, but has an entirely analogous flavour. As before, we use Theorem~\ref{vir pull-back} and Example~\ref{excess} to get that $i_*[X]^{vir}=4\big[H^2\big]$, with $H\in A_*\big(\PP^3\big)$ the class of a hyperplane. Although more elaborate, the computation in Example~\ref{mainexample} has the advantage that it gives a splitting of the virtual class on the components of~$X$.
\end{Remark}
\begin{Definition} Let $X\xrightarrow{f} Y\xrightarrow{g} Z$ be morphisms of DM stacks and let $E^{\bullet}_{X/Y}$, $E^{\bullet}_{X/Z}$, $E^{\bullet}_{Y/Z}$ be obstruction theories for $f,g\circ f, g$ respectively. We call the obstruction theories \emph{compatible} if there is a morphism of exact triangles:
\[\begin{tikzcd}
f^*E^{\bullet}_{Y/Z}\ar[r,"g^{\#}"]\ar[d] & E^{\bullet}_{X/Z}\ar[r,"f^{\#}"]\ar[d] & E^{\bullet}_{X/Y}\ar[r]\ar[d] & f^*E^{\bullet}_{Y/Z}[1]\ar[d] \\
f^*\mathbb L^{\bullet}_{Y/Z}\ar[r] & \mathbb L^{\bullet}_{X/Z}\ar[r] & \mathbb L^{\bullet}_{X/Y}\ar[r] & f^*\mathbb L^{\bullet}_{Y/Z}[1].
\end{tikzcd}
\]
\end{Definition}
The following is a more general version of Proposition~\ref{easycompat}.
\begin{Theorem}[\cite{vir-pull}]
Let $X\xrightarrow{f} Y\xrightarrow{g} Z$ be morphisms of DM stacks and $(E^{\bullet}_{X/Y}, E^{\bullet}_{X/Z},E^{\bullet}_{Y/Z})$ a~triple of compatible obstruction theories. Then
\[(g\circ f)^{!}_{E^{\bullet}_{X/Z}}=f^{!}_{E^{\bullet}_{X/Y}}\circ g^{!}_{E^{\bullet}_{Y/Z}}.\]
\end{Theorem}

\subsection{Push-forwards of virtual classes}
Recall that, to any proper morphism $f\colon X\to Y$ of DM stacks, we can associate a push-forward operation
\[f_*\colon \ A_*(X)\to A_*(Y),\]
defined on prime cycles by
\[f_*[V]=\begin{cases}f_*[V]=\deg(f|_{V})[f(V)], &\text{if $\dim(f(V))=\dim(V)$},\\
0, &\text{otherwise.}
\end{cases}
\]
If $X$ and $Y$ are varieties of the same dimension, then $f_*[X]=\deg(f)[Y]$. If $X$ and $Y$ are pure-dimensional schemes such that the dimension of $X$ is greater or equal than the dimension of~$Y$, and $Y=\bigcup Y_i$ with $Y_i$ irreducible, then $f_*[X]=\sum c_i[Y_i]$, for some $c_i\in\mathbb Q$. If~$X$ is not irreducible, see~\cite[Chapter 1]{Fulton} for a definition of the fundamental class of~$X$.

This discussion is generally false for virtual classes. Let us look at an easy example.
\begin{Example} Consider $\tilde{\mathbb P}^2=\operatorname{Bl}_p\PP^2$ with two different obstruction theories,
\begin{equation*}
 E_1=\big[\mathcal O(-H)\xrightarrow{0}\Omega_{\tilde{\mathbb P}^2}\big] \qquad\text{and}\qquad E_2=\big[\mathcal O(-E)\xrightarrow{0}\Omega_{\tilde{\mathbb P}^2}\big]
\end{equation*}
with corresponding virtual classes given by $\big[\tilde{\mathbb P}^2\big]^{\vv}_1=H$ and $\big[\tilde{\mathbb P}^2\big]^{\vv}_2=E$. Let $\operatorname{id}\colon \tilde{\mathbb P}^2\to \tilde{\mathbb P}^2$ denote the identity. Then, $(\operatorname{id}_{\tilde{\mathbb P}^2})_*[H]\neq[E]$ which shows that
\[\operatorname{id}_* \big[\tilde{\mathbb P}^2\big]^{\vv}_1\neq \big[\tilde{\mathbb P}^2\big]^{\vv}_2.
\]
\end{Example}
\begin{Definition}\label{defvpf} We say that a proper morphism of DM type $f\colon X\to Y$ between stacks with perfect obstruction theories satisfies the \emph{virtual push-forward property} if
 \[f_*[X]^{\vv}=\sum c_i[Y_i]^{\vv}
 \]
 for some splitting of $[Y]^{\vv}=\sum[Y_i]^{\vv}$ and some $c_i\in\mathbb Q$.

 We say that $f$ satisfies the \emph{strong virtual push-forward} property if moreover all $c_i$'s are equal.
\end{Definition}
\begin{Theorem}[\cite{vir-push}] \label{vpf1} Let $f\colon X\to Y$ be a proper morphism of DM type between algebraic stacks. Assume that~$Y$ is connected, and there is a commutative diagram of perfect obstruction theories
 \[
 \begin{tikzcd}
 f^*E^\bullet_{Y/\Spec(k)}\ar[r,"\phi"]\ar[d] & E^\bullet_{X/\Spec(k)} \ar[d] \\
 f^*\mathbb L^\bullet_{Y/\Spec(k)}\ar[r] & \mathbb L^\bullet_{X/\Spec(k)}
 \end{tikzcd}
 \]
such that the cone of $\phi$ is supported in degrees $[-1,0]$. Then $f$ satisfies the strong virtual push-forward property in homology.
\end{Theorem}
In practice, the cone of $\phi$ is supported in degrees $[-1,0]$ iff
\begin{gather*}h^{-1}(\phi)\colon \ h^{-1}\big(f^*E^\bullet_{Y/\Spec(k)}\big)\to h^{-1}\big(E^\bullet_{X/\Spec(k)}\big)\end{gather*} is injective, or dually iff the map \begin{gather*}h^{1}\big(\phi^\vee\big)\colon \ h^{1}\big(E^{\bullet\vee}_{X/\Spec(k)}\big) \to h^{1}\big(f^*E^{\bullet\vee}_{Y/\Spec(k)}\big)\end{gather*} is surjective.

The condition in Theorem~\ref{vpf1} is not necessary.
\begin{Example}\label{twodiv} Let $X$ be the union of two smooth divisors $D_1$, $D_2$ on a smooth space $V$. Then the obvious virtual class is \begin{gather*}[X]^{\vv}=[D_1]+[D_2]\end{gather*} and we see that the embedding $i_1\colon D_1\to X$ satisfies the virtual push-forward property but we do not necessarily have a surjective morphism between obstruction theories. We do have a morphism between obstruction theories $O_{D_1}(D_1)\to O_{D_1}(D_1+D_2)$ given by multiplication with~$D_2$. This is an isomorphism if $D_1$ does not intersect~$D_2$, but it is not surjective if the intersection is non-empty.
\end{Example}

Let us give a new virtual push-forward theorem which is better than Theorem~\ref{vpf1} in the sense that it covers Example~\ref{twodiv}, but has the disadvantage of a hypothesis which may be hard to check in general.
\begin{Theorem}\label{vpfnew}
 Let $X\stackrel{f}{\hookrightarrow} Y\hookrightarrow V$ be embeddings of DM stacks such that $X$ and $Y$ have obstruction bundles $E_{X/V}$ and $E_{Y/V}$ of the same rank. Consider the complexes ${E^\bullet}_{Y/V}=\big[E^\vee_{Y/V}\to 0\big]$ and $E^{\bullet}_{X/V}=\big[E^\vee_{X/V}\to 0\big]$ supported in $[-1,0]$ and assume that there is a commutative diagram of perfect obstruction theories
 \begin{equation*}
 \begin{tikzcd}
 f^*{E^{\bullet}}_{Y/V}\ar[r,"{\phi^\vee}"]\ar[d] & {E^{\bullet}}_{X/V} \ar[d] \\
 f^*\mathbb L^\bullet_{Y/V}\ar[r] & \mathbb L^\bullet_{X/V} .
 \end{tikzcd}
 \end{equation*}
Let $f^*C_{Y/V}$ denote the fibre product $C_{Y/V}\times_YX$, let \begin{gather*}C=f^*C_{Y/V}\times_{f^*E_{Y/V}}E_{X/V}\end{gather*} and consider the induced diagram
\[
\begin{tikzcd}
C_{X/V}\ar[r,"j"] & C\ar[r]\ar[d] & f^*C_{Y/V}\ar[d] \\
 & E_{X/V}\ar[r,"\phi"] & f^*E_{Y/V}.
 \end{tikzcd}\]
\begin{enumerate}\itemsep=0pt
\item[$1.$] The complex \begin{gather*}F^\bullet=[E_{X/V}\to f^*E_{Y/V}]\end{gather*} concentrated in $[0,1]$ is a perfect dual obstruction theory for $\phi$.
\item[$2.$] If $j_*[C_{X/V}]=N\phi_F^![f^*C_{Y/V}]$ for some $N\in\mathbb Q$, then $f$ satisfies the virtual push-forward property.
\end{enumerate}
\end{Theorem}

\begin{proof} The first claim follows from the triangle of cotangent complexes associated to $E_{X/V}\to f^*E_{Y/V}\to X$.

For the second one denote $[Y]^{\vv}|_ X=0_{f^*E_{Y/V}}^![f^*C_{Y/V}]$. Then, by functoriality of virtual pull-backs, we have that
\[[Y]^{\vv}|_ X=0^!_{E_{X/V}}\phi_F^![f^*C_{Y/V}].
\]
By hypothesis we have $j_*[C_{X/V}]=N\phi_F^![f^*C_{Y/V}]$. This shows that
\begin{equation}\label{YNX}[Y]^{\vv}|_ X=\frac{1}{N}[X]^{\vv}.
\end{equation} In the following we prove that
\begin{equation}\label{YVG}f_*[Y]^{\vv}|_ X=\sum c_i[G_i],
\end{equation} for $G_i$ such that $[Y]^\vv=\sum [G_i]$. By commutativity of virtual pullbacks with push-forwards in the following Cartesian diagram
\[
\begin{tikzcd}
X\ar[r]\ar[d] & Y\ar[d] \\ f^*E_{Y/V}\ar[r] & E_{Y/V}
\end{tikzcd}
\]
it is enough to check that $\phi_*[f^*C_{Y/V}]=\sum d_i [\tilde G_i]$, with $\tilde G_i$ the irreducible components of $C_{Y/V}$. This is now clear because $f^*C_{Y/V}$ and $C_{Y/V}$ are both pure of the same dimension.

The conclusion follows from (\ref{YNX}) and~(\ref{YVG}).
\end{proof}

\begin{Remark} Note that $C_{Y/V}\times_Y X$ could have components of dimension smaller than the dimension of $V$. In this case $j_*[C_{X/V}]\neq N\phi_F^![f^*C_{Y/V}]$. The example of~$f$ we have in mind is an embedding of components. Note that in the above proof $c_i$ can be zero.
\end{Remark}

\subsection{Computations of virtual classes} \label{algorithms}
\subsubsection{Computations of virtual classes via Segre classes}
Given a global embedding $i\colon X\to V$ and an obstruction theory $E_{X/V}$, the computation of the virtual class reduces to Segre classes computations (see Proposition~\ref{prop:zerocone}). Segre classes have been computed in many instances (see, e.g., \cite{A1,A, AB,EJP,Helmer, Helmer2, MQ}). Although existent computations are usually for~$i_*s(X,V)$, ideas can probably be adapted
to compute $[X]^{\vv}$ instead of $i_*[X]^{\vv}$.

\subsubsection{Method 1} Consider the blow-up diagram
\[
\begin{tikzcd}
\tilde{X}\ar[r]\ar[d,"\eta"] & \tilde{V}\ar[d] \\
X\ar[r] & V,
\end{tikzcd}
\]
where $X\hookrightarrow V$ is a closed embedding in a smooth scheme, and $\tilde{V}$ is the blow-up of $V$ along~$X$. Then, given an obstruction bundle $E_{X/V}$, $\eta^*E_{X/V}$ is an obstruction bundle for $\tilde{X}$ in $\tilde{V}$. There is an exact sequence of vector bundles
\[ 0\to \mathcal O_{\tilde{X}}\big(\tilde{X}\big)\to \eta^*E_{X/V}\to Q\to 0
\]
and
\[
[X]^{\vv}=\eta_* \big(c_{\text{top}}(Q)\cap \big[\tilde{X}\big]\big).
\]
This method is used in Examples \ref{segrecones} and \ref{exteriorpoint}.
\begin{Remark} This is the basic idea of Aluffi \cite{A1}. Although simple enough, it might be complicated to implement it when $X$ and $V$ are moduli spaces. Blowing them up might mean loosing modularity, and then Chern classes computations might become more difficult.
\end{Remark}
\begin{Remark}Other ideas have been used successfully to compute Segre classes. For example~\cite{EJP} uses residual schemes. As far as we know, these ideas have not been applied to our context yet.
\end{Remark}

\subsubsection{Computations of virtual classes via deformation}

\subsubsection{Method 2} Given an embedding $X\hookrightarrow V$ and a set of global generators $\{f_1,\dots,f_k\}$ for the ideal $I_{X/V}$, we consider the obstruction bundle defined by restricting $\oplus \OO_V(f_i)$ to $X$. Set $Z_{k+1}=V$ and $i_j\colon Z_j=V(f_j)\cap Z_{j+1} \hookrightarrow Z_{j+1}$. Note that $Z_j$ is a union of divisors and components in $Z_{j+1}$. We use this sequence of embeddings to deform $C_{Z_j/V}$ to $C_{Z_j/C_{Z_{j+1}/V}}$. The virtual class on $X$ is then $[X]^{\vv}=i_1^!(\dots i_k^![V])$.
Each pull-back can be computed as in Example \ref{mainexample}. This is easy, but very restrictive: it applies only to (possibly non-transverse) intersections of divisors.

This method is used in Examples \ref{mainexample} and \ref{nonlimit}. It could also be applied to compute the virtual class in Example \ref{doublepoint}.
\subsubsection{Generalisation}
More generally, assume that the obstruction theory of $X$ in $V$ is the restriction of a vector bundle $E$ on $V$ with a surjection $E^\vee\to I_{X/V}$. Assume furthermore that~$E$ admits a filtration $0\subseteq E_1\subseteq\dots \subseteq E_r=E$, with line bundle quotients $L_k=E_k/E_{k-1}$. By dualising the maps induced by the filtration we get diagrams:
\[\begin{tikzcd}
L_k^\vee \ar[r,hook] & E_k^\vee \ar[r,two heads]\ar[d] & E_{k-1}^\vee \\
& i_k^*I_{X/V}\ar[r] & \OO_{V_k}
\end{tikzcd}
\]
hence (possibly vanishing) cosections $f_k\colon L_k^\vee\to \OO_{V_k}$, where $i_k\colon V_k\hookrightarrow V$ is defined iteratively by the vanishing of $\{f_r,\ldots,f_{k+1}\}$. Note that by this definition we have maps $E_{k-1}^\vee\to i_{k-1}^*I_{X/V}$. From Construction \ref{defvb} and Proposition~\ref{easycompat} we get an embedding of a cone inside a split vector bundle \[C_{X/C_{V_1/\ldots_{/C_{V_{r-1}/V}}}}\subseteq \bigoplus_{k=1}^r L_k.\]

\section{Applications to moduli spaces of stable maps}\label{applications}
In this section we illustrate the theory in the previous sections in more meaningful examples. We define moduli spaces of stable maps and look at their irreducible components. We define Gromov--Witten invariants and we discuss the quantum Lefschetz hyperplane property in genus zero. In genus one the quantum Lefschetz hyperplane property fails \cite{coates2012}. One way to fix this problem is by first separating degenerate contributions. We show how to do this in one simple example.

\subsection{Moduli spaces of stable maps and Gromov--Witten invariants}

In this section we define moduli spaces of stable maps, which are useful compactifications of spaces of curves embedded in a given space. For an application of moduli spaces of stable maps to enumerative geometry see~\cite{k}.
\begin{Definition}[\cite{k}] Let $X$ be a smooth projective variety over $\mathbb C$, $\beta\in H_2(X,\mathbb Z)^+$, $(g,n)\in\mathbb N^2$, where $H_2(X,\mathbb Z)^+$ denotes the semigroup of homology classes of algebraic curves modulo torsion. The moduli space of stable maps $\MM{g}{n}{X}{\beta}$ parametrises $[f\colon (C,x_1,\ldots,x_n)\to X]$ such that:
\begin{enumerate}\itemsep=0pt
\item[(1)] $C$ is a projective, connected, reduced, nodal curve of genus $g$,
\item[(2)] $x_1,\ldots,x_n$ are distinct smooth points of $C$,
\item[(3)] the map $f$ satisfies $f_*[C]=\beta$,
\item[(4)] stability condition: any contracted component of geometric genus 0 has at least 3 special points (i.e. nodes or markings), and any contracted curve of genus 1 has at least 1 special point.
\end{enumerate}
\end{Definition}

\begin{Theorem}
 $\MM{g}{n}{X}{\beta}$ is a proper DM stack~{\rm \cite{k}}. Let $X$ be a homogeneous space, e.g., $X=\PP^N$. The stack $\MM{0}{n}{X}{\beta}$ is smooth and irreducible~{\rm \cite{FP}}.
\end{Theorem}

\subsubsection{Obstruction theory from deformation theory} \label{obsdefth}
Let $\mathfrak M_{g,n}$ be the algebraic stack whose points parametrise pre-stable genus $g$ curves with $n$ marked points (see~\cite{Ber}). A deformation theory argument shows that $\mathfrak M_{g,n}$ is a smooth Artin stack of dimension $3g-3+n$. Let $\rho\colon \MM{g}{n}{X}{\beta}\to \mathfrak M_{g,n}$ denote the morphism which forgets the map. In the following we introduce an obstruction theory for $\rho$.
Consider the morphisms
\[\xymatrix{\mathcal C\ar[r]^f\ar[d]^{\pi} & X,\\
				\overline{\mathcal{M}}_{g,n}(X,\beta)}\]
where $\pi\colon \mathcal C\to \MM{g}{n}{X}{\beta}$ is the universal curve and $f\colon \mathcal{C}\to X$ is the universal map.
A relative dual obstruction theory for $\rho$ is given by the complex
\begin{equation}\label{mapsobs}\E^{\bullet}_{\rho}=R^{\bullet}\pi_* f^*T_X.
\end{equation}
We define the \emph{expected dimension} of $\MM{g}{n}{X}{\beta}$ as
\[
 \operatorname{expdim}\big(\MM{g}{n}{X}{\beta}\big)=(1-g)(\dim(X)-3)-K_X\cdot\beta+n.
\]
\subsubsection{Gromov--Witten invariants}The obstruction theory in~(\ref{mapsobs}) induces a virtual class
\[
\big[\overline{\mathcal{M}}_{g,n}(X,\beta)\big]^{\vv}_{\E^{\bullet}_{\rho}}\in A_{\mathrm{expdim}}\big(\overline{\mathcal{M}}_{g,n}(X,\beta)\big).
\]
Consider the evaluation maps $\mathrm{ev}_i\colon \overline{\mathcal{M}}_{g,n}(X,\beta)\to X$, defined by $\mathrm{ev}_i(C, x_1,\dots,x_n,f)=f(x_i)$.
Consider $\gamma_i\in A^{k_i}(X)$ such that $\sum_i{k_i}=\mathrm{expdim}(\overline{\mathcal{M}}_{g,n}(X,\beta))$. We call
\[
\langle \gamma_1,\ldots,\gamma_n\rangle^X_{g,d}:=\deg\big(( \mathrm{ev}_1^*\gamma_1\cup\dots\cup \mathrm{ev}_n^*\gamma_n)\cap \big[\overline{\mathcal{M}}_{g,n}(X,\beta)\big]^{\vv}_{\E^{\bullet}_{\rho}}\big)\in \mathbb{Q}
\]
a Gromov--Witten invariant with insertions $\gamma_i$.

\subsection{Examples of moduli spaces of stable maps}
\begin{Example}
 We have an isomorphism $\MM{g}{n}{X}{0}\simeq \overline{\mathcal M}_{g,n}\times X$, which shows that $\MM{g}{n}{X}{0}$ has dimension $3g-3+n+\dim(X)$. Note that this is bigger than the expected dimension. The difference is $g\dim(X)$, which is the rank of the obstruction bundle $\mathbb E^\vee\boxtimes T_X$. Here, $\mathbb E$ denotes the Hodge bundle $\pi_*\omega_\pi$. The fibres of the obstruction bundle are
 \[H^1(C,f^*T_X)=H^1(C,\OO_C)\otimes T_X.
 \]
 \end{Example}
\begin{Example}
Let us look at $\overline{{\mathcal M}}_{1,0}\big({\mathbb P}^2,3\big)$ and show it has several irreducible components.
\begin{itemize}\itemsep=0pt
 \item The \emph{main} component is the closure of the open locus where the source curve is smooth irreducible and the map is an embedding. Recall that a smooth cubic has genus one by the degree-genus formula. The main component has the same dimension as the space of cubics in~$\mathbb P^2$, that is 9.
 \item $D^1$ is the closure of the locus where the source curve has two components: one rational curve glued to an elliptic curve and the map contracts the elliptic curve and sends the rational curve to a nodal cubic. To compute the dimension of this locus we have to choose a nodal cubic, a point on it corresponding to the image of the genus one curve, and an elliptic curve with one marked point. Adding these dimensions we get \begin{gather*}10=8+1+1.\end{gather*}

 \item $D^2$ is the closure of the locus where the source curve consists of two rational tails mapping to a line and a conic respectively, and an elliptic bridge which gets contracted. To compute the dimension of this locus we have to choose a conic, a line in the plane, and a point of~$\overline{\mathcal M}_{1,2}$. Adding these up, we get \begin{gather*}9=5+2+2.\end{gather*}
 \end{itemize}
 \end{Example}

\begin{Example}\label{ex:components} More generally, $\MM{1}{n}{\mathbb P^r}{d}$ has several irreducible components: a main component $\MM{1}{n}{\mathbb P^r}{d}^{\mathrm{main}}$, which is the closure of the locus of maps from a smooth elliptic curve, and boundary components \begin{gather*}\mathcal{D}(\mathbb{P}^r, d)^{\lambda}=\MM{1}{k+n_0}{\mathbb P^r}{0}\times_{\mathbb{P}^r}\MM{0}{n_1+1}{\mathbb P^r}{d_1}\times_{\mathbb{P}^r}\cdots \times_{\mathbb{P}^r} \MM{0}{n_k+1}{\mathbb P^r}{d_k}.\end{gather*} Here, $\lambda$ denotes the data $(k;d_1,\dots,d_k; n_0,\dots,n_k) $ such that $\sum d_i=d$ and $\sum\limits_{i=0}^k n_i=n$. We have that $\mathcal{D}(\mathbb{P}^r, d)^{\lambda}$ are either components or closed substacks of the main component. See Vakil's smoothability criterion \cite[Lemma~5.9]{Vre} for $k>r$, or the discussion in \cite[Section~2]{BCM}. In the following we denote by $\mathcal{D}(\mathbb{P}^r, d)^1$ the boundary component $\mathcal{D}(\mathbb{P}^r, d)^{\lambda}$ with $\lambda=(1,d,0,n)$. With this notation we have that \begin{gather*}\MM{1}{n}{\mathbb P^r}{d}=\MM{1}{n}{\mathbb P^r}{d}^{\mathrm{main}}\bigcup_{\lambda} \mathcal{D}(\mathbb{P}^r, d)^{\lambda}.\end{gather*}
\end{Example}

\begin{Example}$\overline{{\mathcal M}}_{g,0}\big({\mathbb P}^1,1\big)$
has an empty ``main component'' for $g>0$, since there is no map of degree 1 from a smooth curve of positive genus to $\mathbb P^1$. Yet, this space is not empty -- it contains maps whose source curve is a smooth rational curve $R$ that is mapped isomorphically to the target glued to $l$ contracted branches of positive genus $g_1,\ldots,g_l$ such that $\sum\limits_{i=1}^l g_i=g$. Note that the restriction of the map to $R$ is an embedding and thus the dual graph of the source curve must be a tree around $R$. The dimension of such loci is $l+\sum(3g_i-3+1)=3g-l$.
\end{Example}
\begin{Example}
Let $\tilde{\mathbb P}^2$ be the blow-up of the projective plane at a point $p$. Let $H$ denote the class of the pullback of a general hyperplane class, and $E$ the class of the exceptional divisor. Then $H_2\big(\tilde{\mathbb P}^2,\mathbb Z\big)=\mathbb Z\langle H,E\rangle$ and the effective cone is $\mathbb N\langle H-E, E\rangle$.

Consider
$\overline{{\mathcal M}}_{0,0}\big(\tilde{\mathbb P}^2,3H\big)$. This space has several components as follows.
\begin{itemize}\itemsep=0pt
\item One component of dimension 8, whose general point consists of a rational plane (nodal) cubic that does not pass through~$p$.

\item Another component of dimension 8 isomorphic to
\begin{gather*}
\overline{{\mathcal M}}_{0,1}\big(\tilde{\mathbb P}^2,3H-2E\big)
\times_{\tilde{\mathbb P}^2}
\overline{{\mathcal M}}_{0,1}\big(\tilde{\mathbb P}^2,2E\big).
\end{gather*} The general point consists of maps from a curve with two components, one mapping in~$\PP^2$ to a cubic with a node in $p$, and the other one covering the exceptional divisor $2:1$. Note that $\overline{{\mathcal M}}_{0,1}\big(\tilde{\mathbb P}^2,2E\big)
\simeq \overline{{\mathcal M}}_{0,1}\big({\mathbb P}^1,2\big) $ because $E$ is rigid. $
\overline{{\mathcal M}}_{0,1}\big({\mathbb P}^1,2\big)$ has dimension~3, but here it comes with a non-trivial obstruction bundle.
\end{itemize}
\end{Example}
\begin{Remark}We conclude this subsection with the following remarks.
\begin{itemize}\itemsep=0pt
 \item For general $g$ and target $X$, the space of stable maps may have many components of different dimensions.
 \item It is difficult to understand what lies in the boundary of the main component and consequently the tangent space (or the intrinsic normal cone) at those points.
\end{itemize}
\end{Remark}

\subsection[Equations for the moduli space of maps to $\mathbb P^r$]{Equations for the moduli space of maps to $\boldsymbol{\mathbb P^r}$} \label{eq from twist}

 \begin{notation} Consider the Artin stack $\mathfrak{M}^{\rm div}_{1,n}$ parametrising genus~1 pre-stable curves with $n$ marked points and a Cartier divisor $D$ of degree $d$ contained in the smooth locus of~$C$. We denote by $\mathfrak{C}_{\mathfrak{M}^{\rm div}}$ the universal curve.
 \end{notation}

We want to find a local embedding of $\MM{1}{n}{\mathbb P^r}{d}$ into a smooth stack~$V$, and equations for $\MM{1}{n}{\mathbb P^r}{d}\hookrightarrow V$. This was done in \cite{CF-K,HL}. We sketch the construction and refer to \cite{HL} for details.

The moduli space $\MM{1}{n}{\mathbb P^r}{d}$ can be constructed as an open subset of the moduli space that parametrises \begin{gather*}(C;x_1,\ldots,x_n;L;s_0,\ldots,s_r),\end{gather*} with $L$ a line bundle on $C$ and $s_i\in H^0(C,L)$. On the open locus $U_0\subseteq \MM{1}{n}{\mathbb P^r}{d}$ where $D=\{s_0=0\}$ does not contain nodes or components, we have a projection
$U_0\to\mathfrak{M}^{\rm div}_{1,n}$. On a chart $\mathcal V$ of $\mathfrak{M}^{\rm div}_{1,n}$ we may pick a section of the universal curve $\mathfrak C_{\mathfrak{M}^{\rm div}}$ that passes through the smooth locus of the minimal subcurve of genus~1. Denote by $\mathcal A$ the corresponding relative Cartier divisor. Let $\mathcal L$ be the universal line bundle on $\mathfrak{C}_{\mathfrak{M}^{\rm div}}$ and consider the exact sequence on~$\mathfrak{C}_{\mathcal V}$:
\[
 0\to \mathcal L\to \mathcal L(\mathcal A)\to \mathcal L(\mathcal A)|_{\mathcal A}\to 0.
\]
By pushing it forward along the universal curve it gives us an exact sequence of sheaves on $\mathcal V$
\begin{equation*}
 0\to \pi_*\mathcal L\to\pi_*\mathcal L(\mathcal A)\to \pi_*\mathcal L(\mathcal A)|_{\mathcal A}\to R^1\pi_*\mathcal L\to 0.
\end{equation*}
Then $\pi_*\mathcal L(\mathcal A)$ is locally free on $\mathcal V$, and we define \begin{gather*}V=\bigoplus _{1}^r\pi_*\mathcal L(\mathcal A).\end{gather*} The exact sequence above can be pulled back to $V$ and we get a tautological section of $\pi_*\mathcal L(\mathcal A)|_{\mathcal A}^{\oplus r}$, whose vanishing locus contains the local chart $U_0$ of $\MM{1}{n}{\mathbb P^r}{d}$ as an open substack.

If we look around a point of the form: a contracted elliptic curve $E$ glued to one rational tail~$R$, it turns out that there are local coordinates $y_i,x_{1,i},\ldots,x_{d,i}$, $i=1,\ldots,r$ on $V$ such that the equations of $U_0\subseteq \MM{1}{n}{\mathbb P^r}{d}$ inside~$V$ are of the form
\begin{equation}\label{genoneeq}
\{\xi y_1=\dots = \xi y_r=0\},
\end{equation}
where $\xi$ is a local coordinate on $\mathcal V$ which vanishes on the locus where the node separating $E$ from $R$ is not smoothed out.

Around a point of the form: contracted elliptic curve with $k$ rational tails of positive degree, equations look more complicated, but after a suitable blow-up, they can be brought to the form above. See \cite{HL} for a complete discussion.
\begin{Remark}The idea of twisting a line bundle on a family of curves with a sufficiently high power of a Cartier divisor appears in \cite{CF-K,CF-K-M}. In this way one gets a global embedding, but equations for this embedding are not easy to find.

By twisting with just one appropriately chosen section one can get a smaller smooth ambient space, together with very explicit equations. This is the approach in \cite[Theorem~4.16]{HL}. The drawback is that this construction is only local on $\MM{1}{n}{\mathbb P^r}{d}$.
\end{Remark}

\subsection{Obstruction theory from local equations} \label{obs:loceq}
Note that the embedding $U_0\hookrightarrow V$, endows $U_0$ with a perfect obstruction theory relative to $\mathfrak{M}^{\rm div}_{1,n}$
\[[\pi_*\mathcal L(\mathcal A)\to \pi_*\mathcal L(\mathcal A)|_{\mathcal A}].\]
We give explicit gluing maps for the above local obstruction theories. We view this construction as an alternative definition of the obstruction theory for $\MM{1}{n}{\mathbb P^r}{d}$.

Fix $\mathcal{A}_1$, $\mathcal{A}_2$, $\mathcal{A}_3$ and $\tilde{\mathcal{A}}$ large enough such that $R^1\pi_*\mathcal{L}(\mathcal{A}_i)=0$ and such that $\tilde{\mathcal{A}}-\mathcal{A}_i$ is effective for every $i$. Denote by $A_i^{\bullet}$ the complex $[\pi_*\mathcal{L}(\mathcal{A}_i)\to \pi_*\mathcal{L}(\mathcal{A}_i)|_{\mathcal{A}_i}]$ concentrated in $[0,1]$ and by $\tilde{A}^{\bullet}$ the complex $\big[\pi_*\mathcal{L}\big(\tilde{\mathcal{A}}\big)\to \pi_*\mathcal{L}\big(\tilde{\mathcal{A}}\big)|_{\tilde{\mathcal{A}}}\big]$ concentrated in $[0,1]$. Consider commutative diagrams
\[\xymatrix{\mathcal{L}\ar[r]\ar[d]&\mathcal{L}\ar[d]\\
\mathcal{L}(\mathcal{A}_i)\ar[r]\ar[d]&\mathcal{L}\big(\tilde{\mathcal{A}}\big)\ar[d]\\
\mathcal{L}(\mathcal{A}_i)|_{\mathcal{A}_i}\ar[r]&\mathcal{L}\big(\tilde{\mathcal{A}}\big)|_{\tilde{\mathcal{A}}}.
}
\]
Take push-forwards via $\pi$ of the morphisms above. The snake lemma, applied to the induced diagram, implies that $A_i^{\bullet}\to \tilde{A}^{\bullet}$ are quasi-isomorphisms. These quasi-isomorphisms induce quasi-isomorphisms $f_{ij}\colon A_i^{\bullet}\to A_j^{\bullet}$ via the roof
\[ \xymatrix{
&&A_j^{\bullet}\ar[ld]\ar[rd]\\
A_i^{\bullet}\ar[r]&\tilde{A}^{\bullet}&&A_j^{\bullet}
}
\]
and it is easy to see that these quasi-isomorphisms satisfy the cocycle condition. Then $\bigoplus^r A_i^{\bullet}$ together with the gluing morphisms $f_{ij}$ define a perfect obstruction theory for $\MM{1}{n}{\mathbb P^r}{d}$ relative to $\mathfrak{M}^{\rm div}_{1,n}.$

\subsection{Genus 0 invariants of a hypersurface in projective space.}

Let $s\in H^0(\mathbb P^r,\mathcal O_{\mathbb P^r}(a))$ be a generic homogeneous polynomial of degree $a$, and let $X=Z(s)$ be the hypersurface in $\PP^r$ cut out by $s$. Note that the inclusion $i\colon X\hookrightarrow\mathbb P^r$ induces a morphism $j\colon \MM{g}{n}{X}{\beta}\hookrightarrow \MM{g}{n}{\mathbb P^r}{i_*\beta}$. We are interested in computing $j_*[\MM{g}{n}{X}{\beta}]^\text{vir}$. It turns out that the answer is simple when $g=0$ and particularly hard in higher genus.

When $g=0$, $\mathcal{E}=\pi_*f^*\mathcal O_{\mathbb P^r}(a)$ is a vector bundle on $\MM{0}{n}{\mathbb P^r}{d}$. Furthermore, $s$ induces a~section $\eta$ of $\mathcal{E}$. In each point $(C,u_0,\ldots,u_r)$, with $(u_0,\ldots,u_r)\in H^0( C,f^*\mathcal O(1))$, $\eta$ is $s(u_0,\ldots,u_r)$. This shows that $\eta$ vanishes precisely on the locus of stable maps to~$\mathbb P^r$ whose image lies inside~$X$.

\begin{Theorem}[\cite{KKP}] Let $\MM{0}{n}{X}{d}=\bigsqcup_{\{\beta\colon i_*\beta=d\}}\MM{0}{n}{X}{\beta}$. Then,
 \[
 j_*\big[\,\MM{0}{n}{X}{d}\big]^{\vv}=c_{ad+1}(\mathcal{E})\cap\big[\,\MM{0}{n}{\mathbb P^r}{d}\big].
 \]
\end{Theorem}
This can be seen as a functoriality statement in intersection theory in the following way. We have a~commutative diagram
\[
\begin{tikzcd}
 \MM{0}{n}{X}{d} \ar[rr,"j"]\ar[dr,"\rho_X" below] & & \MM{0}{n}{\mathbb P^r}{d} \ar[dl,"\rho_\mathbb P"] \\
 & \mathfrak{M}_{0,n} &
\end{tikzcd}
\]
equipped with a compatible triple of dual obstruction theories as follows. The exact sequence of vector bundles on $X$
\[
 0\to T_X\to T_{\mathbb P^r|X}\to N_{X/\mathbb P^r}\simeq \mathcal O_X(a)\to 0
\]
induces the following exact triangle on $\MM{0}{n}{X}{d}$
\begin{equation}\label{bigcompat}
 \begin{tikzcd}
R^\bullet\pi_*f^* (T_X)\ar[r] & R^\bullet\pi_*f^* (T_{\mathbb P^r|X})\ar[r] & R^\bullet\pi_*f^* (N_{X/\mathbb P^r})\ar[r] & R^\bullet\pi_*f^* (T_X)[1]\\
\mathbb T_{\overline{\mathcal M}(X)/\mathfrak M}\ar[r]\ar[u] & j^*\mathbb T_{\overline{\mathcal M}(\mathbb P^r)/\mathfrak M}\ar[r]\ar[u] & \mathbb T_{\overline{\mathcal M}(X)/\overline{\mathcal M}(\mathbb P^r)}[1]\ar[r]\ar[u] & \mathbb T_{\overline{\mathcal M}(X)/\mathfrak M}[1]\ar[u].
 \end{tikzcd}
\end{equation}
Here we have denoted $\MM{0}{n}{X}{d}$ by $\overline{\mathcal M}(X)$, $\MM{0}{n}{\mathbb P^r}{d}$ by $\overline{\mathcal M}(\mathbb P^r)$ and $\mathfrak M_{0,n}$ by $\mathfrak M$. Moreover, we have
\begin{enumerate}\itemsep=0pt
 \item $R^\bullet\pi_*f^* (T_{\mathbb P^r|X})\simeq j^*R^\bullet\pi^{\mathbb P}_*f_{\mathbb P}^* (T_{\mathbb P^r})$ (it follows from flat base-change in the derived category); this complex is supported in degree 0.
 \item $R^\bullet\pi_*f^* (N_{X/\mathbb P^r})\simeq \mathcal{E}$ is a vector bundle, supported in degree 0, and it gives a dual perfect obstruction theory for $j$.
\end{enumerate}
Let $\rho_{\PP}^!$ and $\rho_{X}^!$ be the virtual pullbacks with respect to the compatible obstruction theories in~(\ref{bigcompat}), and let $j^!$ be the Gysin pullback induced by the section $\eta\in H^0\big(\MM{0}{n}{\mathbb P^r}{d},\mathcal{E}\big)$. By~(1) we have
\[\big[\,\MM{0}{n}{\mathbb P^r}{d}\big]^{\vv}=\big[\,\MM{0}{n}{\mathbb P^r}{d}\big].
\]
By functoriality of virtual pullbacks~\cite{vir-pull} we have
\[
 \big[\,\MM{0}{n}{X}{d}\big]^\text{vir}=\rho_X^![\mathfrak M_{0,n}]=j^!\rho_{\mathbb P}^![\mathfrak M_{0,n}]=j^!\big[\,\MM{0}{n}{\mathbb P^r}{d}\big].
\]
By Theorem~\ref{vir pull-back} and Example~\ref{excess} applied to the Cartesian diagram
\[\xymatrix{\MM{0}{n}{X}{d}\ar[r]\ar[d]_j&\MM{0}{n}{\PP}{d}\ar[d]^s\\
\MM{0}{n}{\PP}{d}\ar[r]^-0 &\mathcal{E}
}
\]
we get $j_*\big[\,\MM{0}{n}{X}{d}\big]^{\vv}=c_{ad+1}(\mathcal{E})\cap\big[\,\MM{0}{n}{\mathbb P^r}{d}\big]$.

\subsection{Genus 1 invariants of a hypersurface in projective space}
The above story fails in higher genus. From now on we look at $g=1$. Unlike for $g=0$, the sheaf $\mathcal{E}=\pi_*f^*\mathcal O_{\mathbb P^r}(a)$ for $a\geq 1$ is not a vector bundle. $\mathcal{E}$ has rank $ad$ on the open consisting of maps with smooth domain and it has rank $ad+1$ on $D(\PP^r,d)^\lambda$. To see this, we look at $H^1(C,f^*T_{\mathbb P^r})$. On a smooth genus one curve, this group is zero. If $E$ is a contracted curve of genus one, then $f^*\mathcal O(a)|_E\simeq \mathcal O_E$ is trivial and $H^1(E,\mathcal O_E)\simeq \mathbb C$. In terms of obstruction theories, this means that the dual relative obstruction theory of \begin{gather*}j\colon \ \MM{1}{n}{X}{d}\to \MM{1}{n}{\mathbb P^r}{d}\end{gather*} is not supported in degree $1$, but in degrees~$[1,2]$.

\subsubsection{The Vakil--Zinger desingularisation}\label{vzdesing} A first approach to define invariants with a better enumerative behaviour was proposed by Zinger, Li--Zinger and Vakil--Zinger \cite{LZ,VZ,zingerstvsred,redgone,zingred}. The rough idea is to blow-up the boundary components until we find a desingularisation of the main component. Moreover, the object defined analogously to $\mathcal{E}$ on this desingularisation is a vector bundle. For this we first introduce an additional stack of decorated pre-stable curves.

\begin{construction}\label{lboundary}Let us define the moduli space of weighted pre-stable curves $\mathfrak M_{1,n}^\text{wt}$: this is an Artin stack which parametrises the data $(C,x_1,\dots,x_n, \underline{d})$, where $C$ is a pre-stable curve, $x_1,\dots,x_n$ markings, and if $C$ has $j$ irreducible components $\underline{d}=(d_1,\dots,d_j)\in\mathbb{Z}^j$. Then $\mathfrak M_{1,n}^\text{wt}$ is \'{e}tale (but not separated) over $\mathfrak M_{1,n}$.

Again, for any partition $\lambda$ of $(d=\sum\underline{d},n)$, we can define $\mathfrak D^\lambda$ the closed substack of $\mathfrak M_{1,n}^\text{wt}$ the general point of which is a smooth elliptic curve of weight~0, with~$k$ rational tails attached to it, with weight and markings distributed according to~$\lambda$. In fact, this is finer than necessary, and we may as well group all the $\lambda$ with the same number of positive rational tails~$k$ together into~$\mathfrak D^k$. Note that these are closed substacks of $\mathfrak M_{1,n}^\text{wt}$ of codimension~$k$.
\end{construction}

We now briefly describe the Vakil--Zinger desingularisation following \cite{VZ} and \cite{HL}. The idea is to start with the moduli space of pre-stable weighted curves $\mathfrak M_{1,n}^\text{wt}$ and construct an Artin stack~$\widetilde{\mathfrak M}_{1,n}^\text{wt}$ by an iterated blow-up. We first blow up $\mathfrak D^1$ inside $\mathfrak M_{1,n}^\text{wt}$ (this has no effect), then the strict transform of $\mathfrak D^2$, which is smooth, and so on. Note that the total weight is constant on connected components of $\mathfrak M_{1,n}^\text{wt}$, and for a fixed weight it is necessary to perform only finitely many blow-up steps. The Vakil--Zinger space is then constructed by taking the following Cartesian diagram
\[
 \begin{tikzcd}
 \widetilde{\mathcal M}_{1,n}(\mathbb P^r,d) \ar[r]\ar[d] & \MM{1}{n}{\mathbb P^r}{d} \ar[d] & f\colon \ (C,x)\to\mathbb P^r \ar[d,mapsto] \\
 \widetilde{\mathfrak M}_{1,n}^\text{wt}\ar[r] & \mathfrak M_{1,n}^\text{wt}, & (C,x,\deg(f^*\mathcal O(1))).
 \end{tikzcd}
\]
Here, the right vertical map is obtained by forgetting the map, but remembering the degree of~$f^*\mathcal O_{\mathbb P^r}(1)$ on the irreducible components of $C$.
Let us introduce more notation. Let
\[\widetilde{\mathcal D}_{1,n}(\mathbb P^r,d)^{\lambda}= \mathcal D_{1,n}(\mathbb P^r,d)^{\lambda}\times _{ \mathcal M_{1,n}(\mathbb P^r,d)} \widetilde{\mathcal M}_{1,n}(\mathbb P^r,d),
\]
with $\mathcal D_{1,n}(\mathbb P^r,d)^{\lambda}$ as in \ref{ex:components}. We denote by $ \widetilde{\mathcal D}_{1,n}(\mathbb P^r,d)^1$, the component with $\lambda\vdash(d,n)$ of the form $(1; 0,d; 0,n)$ (See Example~\ref{ex:components} for the notation). Notice that $ \widetilde{\mathcal D}_{1,n}(\mathbb P^r,d)^1$ is birational to $ \mathcal D_{1,n}(\mathbb P^r,d)^1$. For any $X\hookrightarrow \PP^r$, let
\[ \widetilde{\mathcal M}_{1,n}(X,\beta)= \mathcal M_{1,n}(X, \beta)\times_{\mathcal M_{1,n}(\PP^r, d)} \widetilde{\mathcal M}_{1,n}(\mathbb P^r,d)\] and
\[\widetilde{\mathcal D}_{1,n}(X,d)^{\lambda}= \mathcal M_{1,n}(X,d)\times _{ \mathcal M_{1,n}(\mathbb P^r,d)} \widetilde{\mathcal D}_{1,n}(\mathbb P^r,d)^{\lambda}.
\] Similarly, we denote $\widetilde{\mathcal D}_{1,n}(X,\beta)^1$, the component with $\lambda\vdash(d,n)$ of the form $(1; 0,d; 0,n)$.
\begin{Theorem}\quad\begin{enumerate}\itemsep=0pt
\item[$(i)$] The blow-up of the main component $\widetilde{\mathcal {M}}_{1,n}(\mathbb P^r,d)^{\mathrm{main}}$ is smooth.
\item[$(ii)$] Let $\pi^{\circ}\colon \mathcal{C}^{\mathrm{main}}\to\widetilde{\mathcal {M}}_{1,n}(\mathbb P^r,d)^{\mathrm{main}}$ be the universal curve over the main component of $\widetilde{\mathcal {M}}_{1,n}(\mathbb P^r,d)$ and let $f^{\circ}\colon \mathcal{C}^{\mathrm{main}}\to \mathbb P^r$ be the universal morphism. Then $\mathcal{E}_{ad}:=\pi^\circ_* (f^\circ)^*\mathcal O_{\mathbb P^r}(a)$ is a vector bundle of rank~$ad$.
\end{enumerate}
\end{Theorem}

For a hypersurface $X\subseteq \mathbb P^r$ of degree $a$, the main component is defined by the Cartesian diagram
\[
 \begin{tikzcd}
 \MM{1}{n}{X}{d}^\text{main} \ar[r]\ar[d] & \MM{1}{n}{\mathbb P^r}{d}^\text{main} \ar[d] \\
 \MM{1}{n}{X}{d} \ar[r] & \MM{1}{n}{\mathbb P^r}{d}
 \end{tikzcd}
\]
and similarly for the Vakil--Zinger desingularisation. Note that $\mathcal{E}_{ad}$ gives a relative perfect obstruction theory for $\tilde{j}\colon \widetilde{\mathcal M}_{1,n}(X,d)^\text{main} \hookrightarrow\widetilde{\mathcal M}_{1,n}(\mathbb P^r,d)^\text{main} $, hence a virtual class
\begin{equation}\label{defreduced}
\tilde{j}_*\big[\widetilde{\mathcal M}_{1,n}(X,d)^\text{main}\big]^{\vv}=c_{ad}(\mathcal{E}_{ad})\cap\big[\widetilde{\mathcal M}_{1,n}({\mathbb P^r},{d})^\text{main}\big]
\end{equation}
of homological degree $d(r+1-a)+n$, where we denoted with~$j$ the embedding into the resolution of the main component of maps to projective space.

\begin{Definition} Let $X\subseteq \mathbb P^r$ be a smooth hypersurface, and $\gamma_1,\ldots,\gamma_n\in A^*(X)$. We call
\[ \langle\gamma_1,\ldots,\gamma_n\rangle_{1,d}^{X,{\rm red}}:=
\deg \big(\big[\widetilde{\mathcal M}_{1,n}(X,d)^\text{main}\big]^{\vv}\cap( \widetilde{\mathrm{ev}}_1^*\gamma_1\cup\dots\cup \widetilde{\mathrm{ev}}_n^*\gamma_n)\big)
\]
the reduced genus~1 invariants of $X$.
\end{Definition}

\begin{Remark} Formula (\ref{defreduced}) can be thought of as a hyperplane property for reduced invariants. Reduced invariants of a Calabi--Yau hypersurface are computable via torus localisation (see \cite{zingred}).
\end{Remark}

\subsubsection{The Zinger/Li--Zinger formula}
The following formula has been proved by Zinger/Li and Zinger~\cite{LZ, zingerstvsred,redgone} in the symplectic category and by Chang and Li~\cite{CL-formula} in the algebraic one.

\begin{Theorem}\label{lizinger}
 Let $X\hookrightarrow\PP^4$ be a smooth hypersurface of degree~$a$, and let $(\gamma_1,\ldots, \gamma_n)\in A^*(X)^{\oplus n}$. Then we have
 \[
 \langle\gamma_1,\ldots,\gamma_n\rangle^X_{1,d}= \langle\gamma_1,\ldots,\gamma_n\rangle_{1,d}^{X,{\rm red}}+\frac{2+d(a-5)}{24} \langle\gamma_1,\ldots,\gamma_n\rangle^X_{0,d}.
 \]
\end{Theorem}

Reasons to believe this formula:
\begin{itemize}\itemsep=0pt
 \item Reduced invariants in genus 1 capture the contribution of the main component to ordinary GW invariants.
 \item From $\mathcal{D}(X,d)^1\simeq \overline{\mathcal M}_{1,1}\times\MM{0}{n+1}{X}{d}$, we get a projection
 \[q^1\colon \ \mathcal{D}(X,d)^1\to \MM{0}{n}{X}{d},
 \]
 which forgets the elliptic curve and the point to which it is attached. The spaces $\MM{0}{n}{X}{d}$ and $\mathcal{D}(X,d)^1$ have the same virtual dimension, so one could hope for some form of virtual push-forward Definition~\ref{defvpf} to hold.

 \item The coefficient can be computed in an ideal situation: the generic fiber is $\overline{\mathcal M}_{1,1}\times R$, where $R\simeq\mathbb P^1$ is the fiber of the universal curve over $\MM{0}{n}{X}{d}$. The obstruction bundle on $\overline{\mathcal M}_{1,1}\times R$ is $\mathbb E^\vee\boxtimes N_{R/X}$. Then the second Chern class gives the coefficient of $GW_0$ in the formula.
 \item The other components $\mathcal{D}(X,d)^{\lambda}$, with $\lambda\neq 1$ should not contribute by (virtual) dimensional reasons.
 \end{itemize}

Reasons why the formula is hard to prove:
\begin{itemize}\itemsep=0pt
 \item As discussed in Section \ref{normalcones} the cone of $\MM{1}{n}{X}{d}$ is not the union of cones of the components of $\MM{1}{n}{X}{d}$. Consequently it is hard to split the virtual class on the components of the moduli space $\MM{1}{n}{X}{d}$.
 \item We do not have a perfect obstruction theory for the boundary components: Let $q$ denote the node which joins the elliptic curve~$E$ and the rational tail~$R$. The vector bundle $\mathbb E^\vee\boxtimes\operatorname{ev}_q^* N_{R/X}$ does not extend to a vector bundle on $\widetilde{\mathcal D}(X,d)^1$ when $q$ is a ramification point for the map.

\item Since the space $\widetilde{\mathcal M}_{1,n}(X,d)$ is rather mysterious, it is hard to deal with the components of its normal cone above. In particular, the normal cone may have components supported at the intersection of the main component and boundary components.

\item Although $\mathcal D(X,d)^1$ fits into the following Cartesian diagram
\begin{equation*}
 \begin{xymatrix}
 {
 \mathcal D(X,d)^1 \ar[r]\ar[d]& \mathcal D\big(\mathbb P^4,d\big)^1\ar[d]^{q^1}\\
 \MM{0}{n}{X}{d} \ar[r] & \overline{\mathcal{M}}_{0,n}\big({\mathbb P}^4,d\big)
 }
 \end{xymatrix}
\end{equation*}
we cannot apply the virtual push-forward theorem. This is because we do not know if $\mathfrak C_{\MM{1}{n}{X}{d}/\mathfrak{M}_{1,n}}|_{\mathcal {D}(X,d)^1}$ maps to $\mathfrak C_{\MM{0}{n}{X}{d}/\mathfrak{M}_{0,n}}$.
\end{itemize}

\begin{Example} Let us look at the components of the virtual class of $\widetilde{\mathcal{M}}_{1,n}(\PP^r,d)$. In the following we show that we have a splitting of the virtual class
\[\big[\widetilde{\mathcal{M}}_{1,n}(\PP^r,d)\big]^{\vv}=\big[\widetilde{\mathcal{M}}_{1,n} (\PP^r,d)^{\mathrm{main}}\big]+\sum_{\lambda}\big[\widetilde{\mathcal{D}}(\PP^r,d)^{\lambda}\big]^{\vv},
\]
where the virtual class on the boundary of the moduli space is given by the restriction of the normal cone $\mathfrak{C}_{\widetilde{\mathcal{M}}_{1,n}(\PP^r,d)/\widetilde{\mathfrak{M}}_{1,n}}$ inside $h^1/h^0(R\pi f^*T_{\PP^r})$ to the boundary components. By the local equations (\ref{genoneeq}), the cone $\mathfrak{C}_{\widetilde{\mathcal{M}}_{1,n}(\PP^r,d)/\widetilde{\mathfrak{M}}_{1,n}}$ does not have embedded components supported on the intersection of the main component with the boundary components. Let us denote the component of $\mathfrak{C}_{\widetilde{\mathcal{M}}_{1,n}(\PP^r,d)/\widetilde{\mathfrak{M}}_{1,n}}$ supported on the main component by $\mathfrak{C}^{\mathrm{main}}_{\widetilde{\mathcal{M}}_{1,n}(\PP^r,d)/\widetilde{\mathfrak{M}}_{1,n}}$ and similarly the component of $h^1/h^0(R\pi f^*T_{\PP^r})$ supported on the main component by $h^1/h^0(R\pi f^*T_{\PP^r})^{\mathrm{main}}$. We need to show that
\[\mathfrak{C}^{\mathrm{main}}_{\widetilde{\mathcal{M}}_{1,n}(\PP^r,d)/\widetilde{\mathfrak{M}}_{1,n}}\hookrightarrow h^1/h^0(R\pi f^*T_{\PP^r})^{\mathrm{main}}\]
 and
 \[\mathfrak{C}_{\widetilde{\mathcal{M}}_{1,n}(\PP^r,d)^{\mathrm{main}}/\widetilde{\mathfrak{M}}_{1,n}}\hookrightarrow h^1/h^0\big(T_{\widetilde{\mathcal{M}}_{1,n}(\PP^r,d)^{\mathrm{main}}/\widetilde{\mathfrak{M}}_{1,n}}\big)
 \]
 give the same class. Let us look at the local picture. As in \ref{eq from twist} we have local embeddings of $\widetilde{\mathcal{M}}_{1,n}(\PP^r,d)$ and $\widetilde{\mathcal{M}}_{1,n}(\PP^r,d)^{\mathrm{main}}$ in a smooth stack $\widetilde{V}$ and local obstruction theories $E_{\widetilde{\mathcal{M}}_{1,n}(\PP^r,d)/\widetilde{V}}$ and $E_{\widetilde{\mathcal{M}}_{1,n}(\PP^r,d)^{\mathrm{main}}/\widetilde{V}}$ respectively. By the local equations of $\widetilde{M}_{1,n}(\PP^r,d)$, we have that
\begin{gather*}
N_{\widetilde{\mathcal{M}}_{1,n}(\PP^r,d)^{\mathrm{main}}/\widetilde{V}}\simeq E_{\widetilde{\mathcal{M}}_{1,n}(\PP^r,d)^{\mathrm{main}}/\widetilde{V}},
\\
C_{\widetilde{\mathcal{M}}_{1,n}(\PP^r,d)/V}^{\mathrm{main}}\simeq E_{\widetilde{\mathcal{M}}_{1,n}(\PP^r,d)/\widetilde{V}}|_{\widetilde{\mathcal{M}}_{1,n} (\PP^r,d)^{\mathrm{main}}}.
\end{gather*}
This shows that the two virtual classes are equal to $[\widetilde{\mathcal{M}}_{1,n}(\PP^r,d)^{\mathrm{main}}]$, so in particular they agree.

One could also use Theorem \ref{vpfnew}, although it is somewhat overcomplicated. From the local the local equations~\eqref{genoneeq} we have local morphisms
\[ E_{\widetilde{\mathcal{M}}_{1,n}(\PP^r,d)^{\mathrm{main}}/\widetilde{V}}\to E_{\widetilde{\mathcal{M}}_{1,n}(\PP^r,d)/\widetilde{V}},
\]
which are multiplication by $\xi$. As in Section~\ref{obs:loceq} these give a global morphism of obstruction theories
\[
E_{\widetilde{\mathcal{M}}_{1,n}(\PP^r,d)^{\mathrm{main}}/\widetilde{\mathfrak{M}}_{1,n}}\to E_{\widetilde{\mathcal{M}}_{1,n}(\PP^r,d)/\widetilde{\mathfrak{M}}_{1,n}}|_{\widetilde{\mathcal{M}}_{1,n} (\PP^r,d)^{\mathrm{main}}},
\]
which maps $\mathfrak{C}^{\mathrm{main}}$ to $\mathfrak{C}_{\widetilde{\mathcal{M}}_{1,n}(\PP^r,d)/\widetilde{\mathfrak{M}}_{1,n}}$. The local discussion above shows that we are under the hypothesis of Theorem \ref{vpfnew}, which gives the result.
\end{Example}

\begin{Example} Let us look at local equations of $\widetilde{\mathcal{M}}_{1,n}(X,d)$ inside $\widetilde{V}$. The easiest case is when $X\simeq \mathbb P^3\subseteq \mathbb P^4$ is a hyperplane. Without loss of generality we may assume that $\mathbb P^3$ is given by $\{z_1=0\}$ in $\PP^4$. Then
\begin{equation}\label{depeq}
 \begin{tikzcd}[cramped]
 U_{\widetilde{\mathcal{M}}_{1,n}(X,d)} \ar[r,hook,"\{y_1=x_{1,1}=\dots=x_{d,1}=0\}" below=.4cm] & U_{\widetilde{\mathcal{M}}_{1,n}(\PP^r,d)} \ar[r,hook," \{\xi y_1=\dots=\xi y_4=0\} " below=.4cm] & \widetilde{V}.
 \end{tikzcd}
\end{equation}

Note that there are $d+1$ equations on the left, and~5 on the right, but $\widetilde{\mathcal{M}}_{1,n}(X,d)$ only has virtual codimension $d+5$ in $\widetilde{V}$. This is because the equation $\xi y_1$ on the right becomes superfluous after adjoining all the equations on the left.
\end{Example}

\begin{Remark}For $X=\PP^3$ as before, let us forget that it is easy to compute the contribution of the boundary of $\widetilde{\mathcal{M}}_{1,n}\big(\PP^3,d\big)$ and try to use the geometry of $\widetilde{\mathcal{M}}_{1,n}\big(\PP^4,d\big)$ to sketch a proof.

We first introduce more notation. Let $\pi\colon \mathcal{C}_0\to \overline{\mathcal{M}}_{0,n}(\PP^3,d)$ be the universal curve and $E_{d+1}$ the vector bundle $\pi_*f^*\mathcal{O}_{\PP^3}(1)$ on $\overline{\mathcal{M}}_{0,n}(\PP^3,d)$. Let $\widetilde{\mathfrak{C}}_{\PP^3}$ be the relative normal cone $\mathfrak{C}_{\widetilde{\mathcal{M}}_{1,n}(\PP^3,d)/\widetilde{\mathfrak{M}}_{1,n}}$. By construction we have an embedding of cone stacks $\widetilde{\mathfrak{C}}_{\PP^3}\hookrightarrow h^1/h^0\big(\widetilde{\pi}_*\widetilde{f}^*T_{\PP^3}\big)$ and therefore $\widetilde{\mathfrak{C}}_{\PP^3}|_{\widetilde{\mathcal{D}}(\PP^3,d)^{\lambda}}\hookrightarrow \big(h^1/h^0(\widetilde{\pi}_*\widetilde{f}^*T_{\PP^3}\big)\big)|_{\widetilde{\mathcal{D}} (\PP^3,d)^{\lambda}}$.

The long exact sequence
\begin{align*}0& \to R^0\widetilde{\pi}_*\widetilde{f}^*T_{\PP^3}\to \tilde{j}^*R^0\widetilde{\pi}_*\widetilde{f}^*T_{\PP^4}\to R^0\widetilde{\pi}_*\widetilde{f}^*N_{\PP^3/\PP^4}\to \\
& \to R^1\widetilde{\pi}_*\widetilde{f}^*T_{\PP^3}\to \tilde{j}^* R^1\widetilde{\pi}_*\widetilde{f}^*T_{\PP^4}\to R^1\widetilde{\pi}_*\widetilde{f}^*N_{\PP^3/\PP^4}\to 0
\end{align*}
restricted to $\widetilde{\mathcal{D}}\big(\PP^3,d\big)^{\lambda}$ is an exact sequence of vector bundles. Here, by abuse of notation we also indicate by $\widetilde{f}$ and $\widetilde{\pi}$ the corresponding maps from the universal curve over $\widetilde{\mathcal{D}}\big(\PP^3,d\big)^{\lambda}$. From the exact sequence
\begin{align*}0\to R^0\widetilde{\pi}_*\widetilde{f}^*T_{\PP^3}\to \tilde{j}^*R^0\widetilde{\pi}_*\widetilde{f}^*T_{\PP^4}\to R^0\widetilde{\pi}_*\widetilde{f}^*N_{\PP^3/\PP^4}
\to R^1\widetilde{\pi}_*\widetilde{f}^*T_{\PP^3}\to R^1\widetilde{\pi}_*\widetilde{f}^*T_{\PP^3}\to 0
\end{align*}
on $\widetilde{\mathcal{D}}\big(\PP^3,d\big)^{\lambda}$, we get a distinguished triangle
\[ \begin{xymatrix}{0\ar[r]\ar[d]&
R^{0}\widetilde{\pi}_*\widetilde{f}^*T_{\PP^3}\ar[r]\ar[d]& \tilde{j}^*R^{0}\widetilde{\pi}_*\widetilde{f}^*T_{\PP^4}\ar[r]\ar[d]&R^{0}\widetilde{\pi}_*\widetilde{f}^*N_{\PP^3/\PP^4}\ar[d]\\
R^{0}\widetilde{\pi}_*\widetilde{f}^*N_{\PP^3/\PP^4}\ar[r]&R^{1}\widetilde{\pi}_*\widetilde{f}^*T_{\PP^3}\ar[r]&R^{1}\widetilde{\pi}_*\widetilde{f}^*T_{\PP^3}\ar[r]&0
.}
\end{xymatrix}
\]
The complexes in the above diagram are written vertically and supported in degrees $[0,1]$. By \cite[Proposition~2.7]{BerFan}, we get an exact sequence of cones
\[R^{0}\widetilde{\pi}_*\widetilde{f}^*N_{\PP^3/\PP^4}\to
h^1/h^0\big(\widetilde{\pi}_*\widetilde{f}^*T_{\PP^3}\big)\to \left[\frac{h^1\pi_*f^*T_{\PP^3}}{\tilde{j}^*h^0\widetilde{\pi}_*\widetilde{f}^*T_{\PP^4}}\right]
\]
and thus a deformation of $h^1/h^0\big(\widetilde{\pi}_*\widetilde{f}^*T_{\PP^3}\big)|_{\widetilde{\mathcal{D}}(\PP^3,d)^{\lambda}}$ to
\begin{equation}\label{limvbstack}\widetilde{\pi}_*\tilde{f}^*N_{\PP^3/\PP^4}\oplus \left[\frac{h^1\pi_*f^*T_{\PP^3}}{\tilde{j}^*h^0\widetilde{\pi}_*\widetilde{f}^*T_{\PP^4}}\right]
\end{equation}
on the boundary of $\widetilde{\mathcal{M}}_{1,n}\big(\PP^3,d\big)$. By~\cite{HL}, we have that
\[
h^1/h^0\big(\widetilde{\pi}_*\widetilde{f}^*T_{\PP^3}\big)|_{\widetilde{\mathcal{D}}(\PP^3,d)^{\lambda}}\simeq \big(h^1/h^0\big(\widetilde{\pi}_*\widetilde{f}^*T_{\PP^3}\big)\big)|_{\widetilde{\mathcal{D}} (\PP^3,d)^{\lambda}}.
\]
In the formula above, $\widetilde{f}$ and $\widetilde{\pi}$ in the first term are maps from the universal curve over $\widetilde{\mathcal{D}}\big(\PP^3,d\big)^{\lambda}$, and in the second, they are maps from the universal curve over $\widetilde{\mathcal{M}}\big(\PP^3,d\big)$. This shows that we have an embedding of cone $\widetilde{\mathfrak{C}}_{\PP^3}|_{\widetilde{\mathcal{D}}(\PP^3,d)^{\lambda}}\hookrightarrow h^1/h^0\big(\widetilde{\pi}_*\widetilde{f}^*T_{\PP^3}\big)|_{\widetilde{\mathcal{D}}(\PP^3,d)^{\lambda}}$. Let $\widetilde{\mathfrak{C}}_{\PP^3}^{\mathrm{lim}}|_{\widetilde{\mathcal{D}}(\PP^3,d)^{\lambda}}$ be the limit of $\widetilde{\mathfrak{C}}_{\PP^3}|_{\widetilde{\mathcal{D}}(\PP^3,d)^{\lambda}}$ in the vector bundle stack~(\ref{limvbstack}).

Let us now look at the local structure of $\widetilde{\mathfrak{C}}_{\PP^3}^{\mathrm{lim}}|_{\widetilde{\mathcal{D}}(\PP^3,d)^{\lambda}}$. Let $\widetilde{U}_0$ and $\widetilde{V}$ analogous to $U_0$ and $V$ in Section~\ref{eq from twist} and let $U_{\widetilde{\mathcal{M}}_{1,n}(\PP^3,d)}=\widetilde{U}_0\cap \widetilde{\mathcal{M}}_{1,n}\big(\PP^3,d\big)$. Let $\widetilde{W}$ in $\widetilde{V}$ be given by equations
\[\langle y_2\xi,y_3\xi,y_4\xi \rangle
\]
and let $\tilde{k}\colon U_{\widetilde{\mathcal{M}}_{1,n}(\PP^3,d)}\to \widetilde{W}$ be the induced local embedding. Using an argument similar to the one in Example~\ref{nonlimit}, we see from equations \ref{depeq} that $\widetilde{\mathfrak{C}}_{\PP^3}^{\mathrm{lim}}|_{\widetilde{\mathcal{D}}(\PP^3,d)^{\lambda}}$ is locally isomorphic to
\[
\big(C_{U_{\widetilde{M}_{1,n}(\PP^3,d)}/\widetilde{W}}\times C_{\widetilde{W}/\widetilde{\mathfrak{M}}_{1,n}}\big)|_{\widetilde{\mathcal{D}}(X,d)^{\lambda}\cap U_{\widetilde{M}_{1,n}(\PP^3,d)}}.
\]
Notice that by construction $\widetilde{W}$ has boundary components isomorphic to $\widetilde{\mathcal{D}}\big(\PP^4,d\big)^{\lambda}\cap \widetilde{V}$, with~$\lambda$ as in Example~\ref{ex:components}. Let $\widetilde{W}^{\lambda}$ denote the boundary components of~$\widetilde{W}$. By the definition of $\widetilde{W}$ we see that we have a Cartesian diagram
\[\xymatrix{
\widetilde{D}_{1,n}\big(\PP^3,d\big)^{\lambda}\cap \widetilde{W}^{\lambda}\ar[r]\ar[d] &\widetilde{W}^{\lambda}\ar[d]^{q^{\lambda}}\\
\overline{\mathcal{M}}_{0,n-n_0}\big(\PP^3,d\big)\ar[r]&\overline{\mathcal{M}}_{0,n-n_0}\big(\PP^4,d\big),
}
\]
with $n_0$ as in Example~\ref{ex:components} and $q^{\lambda}$ the map which contracts elliptic components (see the discussion after Theorem~\ref{lizinger} for~$q^1$). This shows that we have a map
\[C_{\widetilde{D}_{1,n}(\PP^3,d)^{\lambda}\cap \widetilde{W}^{\lambda}/\widetilde{W}^{\lambda}}\to C_{\overline{\mathcal{M}}_{0,n-n_0}(\PP^3,d)/\overline{\mathcal{M}}_{0,n-n_0}(\PP^4,d)}
\]
and therefore a map
\[
\widetilde{\mathfrak{C}}^{\mathrm{lim}}_{\PP^3}|_{\widetilde{\mathcal{D}}(\PP^3,d)^{\lambda}}\to C_{\overline{\mathcal{M}}_{0,n-n_0}(\PP^3,d)/\overline{\mathcal{M}}_{0,n-n_0}(\PP^4,d)}.
\]
Using an argument as in \cite[Lemma~4.4]{vir-push}, this shows that the contribution of $\widetilde{\mathfrak{C}}_{\PP^3}$ on the boundary satisfies the virtual push forward property
\[q^{\lambda}_*\big[\widetilde{\mathcal{D}}\big(\PP^3,d\big)^{\lambda}\big]^{\vv}=N^{\lambda} \big[\,\overline{\mathcal{M}}_{0,n-n_0}\big(\PP^3,d\big)\big],
\]
for some $N^{\lambda}\in \mathbb Q$. The constants~$N^{\lambda}$ can be computed as in the discussion after Theorem~\ref{lizinger}.

The situation is similar for hypersurfaces of higher degree, but there are two main complications: the map $q^{\lambda}$ is not proper anymore and the equations on the left are much less clear. In particular, we cannot know whether they are independent of the set of equations on the right and we cannot find an analogous local description of the limiting cone. These are the main reasons why the above idea does not generalise to higher degree. In~\cite{eu3} T.~Coates and the third-named author substantially modify the above proof.
\end{Remark}

\subsection*{Acknowledgements} Most of the material presented here was presented at the Summer School in Enumerative Geometry at SISSA, Trieste, in 2017. A shorter version was part of a Lecture series within the Workshop on Moduli and Mirror Symmetry held by C.M.\ in Alpensia, in 2016. We thank the organisers of these schools for the invitation and wonderful work environment. We would also like to thank the participants for their interest~-- part of the material in this note grew to answer the questions asked during these workshops. We thank Tom Coates for comments on the manuscript. We also thank the referees for their attentive remarks.

L.B.\ was supported by a Royal Society 1st Year URF and DHF Research Grant Scheme.
C.M.\ is supported by an EPSRC funded Royal Society Dorothy Hodgkin Fellowship.
This work was supported by the Engineering and Physical Sciences Research Council grant EP/L015234/1: the EPSRC Centre for Doctoral Training in Geometry and Number Theory at the Interface. This note was finished while L.B.\ and C.M.\ were in residence at the Mathematical Sciences Research Institute in Berkeley, California, during the Spring 2018 semester and it is supported by the National Science Foundation under Grant No.~DMS-1440140.

\pdfbookmark[1]{References}{ref}
\LastPageEnding

\end{document}